\documentclass[12pt,a4paper]{article}
\usepackage{amssymb}
\usepackage{makeidx}
\usepackage[utf8]{inputenc}
\usepackage[T1]{fontenc}
\usepackage[british]{babel}
\usepackage[round]{natbib}
\usepackage{amsmath}
\usepackage{amsfonts}
\usepackage{amsbsy}
\usepackage{bm}
\usepackage{amsthm}
\usepackage[arrow,matrix,pdf]{xy}

\theoremstyle{plain}
\newtheorem{theorem}{Theorem}
\newtheorem{lemma}{Lemma}
\theoremstyle{definition}
\newtheorem{definition}{Definition}
\newtheorem{example}{Example}

\begin{document}

\title{Mat\'{e}rn Class Tensor-Valued Random Fields and Beyond}
\author{Nikolai Leonenko\thanks{%
Cardiff University, United Kingdom} \and Anatoliy Malyarenko\thanks{%
Mälardalen University, Sweden}}
\date{\today }
\maketitle

\begin{abstract}
We construct classes of homogeneous random fields on a three-dimensional
Euclidean space that take values in linear spaces of tensors of a fixed rank
and are isotropic with respect to a fixed orthogonal representation of the
group of $3\times 3$ orthogonal matrices. The constructed classes depend on
finitely many isotropic spectral densities. We say that such a field belong
to either the Mat\'{e}rn or the dual Mat\'{e}rn class if all of the above
densities are Mat\'{e}rn or dual Mat\'{e}rn. Several examples are considered.
\end{abstract}

\section{Introduction}

Random functions of more than one variable, or \emph{random fields}, were
introduced in the 20th years of the past century as mathematical models of
physical phenomena like turbulence, see, e.g., \citet{Friedmann1924}, %
\citet{Karman1938}, \citet{MR0001702}. To explain how random fields appear
in continuum physics, consider the following example.

\begin{example}
Let $E=E^3$ be a three-dimensional Euclidean point space, and let $V$ be the
translation space of $E$ with an inner product $(\boldsymbol{\cdot},%
\boldsymbol{\cdot})$. Following \cite{MR1162744}, the elements $A$ of $E$
are called the \emph{places} in $E$. The symbol $B-A$ is the vector in $V$
that translates $A$ into $B$.

Let $\mathcal{B}\subset E$ be a subset of $E$ occupied by a material, e.g.,
a turbulent fluid or a deformable body. The temperature is a rank~$0$
tensor-valued function $T\colon\mathcal{B}\to\mathbb{R}^1$. The velocity of
a fluid is a rank~$1$ tensor-valued function $\mathbf{v}\colon\mathcal{B}\to
V$. The strain tensor is a rank~$2$ tensor-valued function $\varepsilon\colon%
\mathcal{B}\to\mathsf{S}^2(V)$, where $\mathsf{S}^2(V)$ is the linear space
of symmetric rank~$2$ tensors over $V$. The piezoelectricity tensor is a
rank~$3$ tensor-valued function $\mathsf{D}\colon\mathcal{B}\to\mathsf{S}%
^2(V)\otimes V$. The elastic modulus is a rank~$4$ tensor-valued function $%
\mathsf{C}\colon\mathcal{B}\to\mathsf{S}^2(\mathsf{S}^2(V))$. Denote the
range of any of the above functions by $\mathsf{V}$. Physicists call $%
\mathsf{V}$ the \emph{constitutive tensor space}. It is a subspace of the
tensor power $V^{\otimes r}$, where $r$ is a nonnegative integer. The form
\begin{equation*}
(\mathbf{x}_1\otimes\cdots\otimes\mathbf{x}_r,\mathbf{y}_1\otimes\cdots%
\otimes\mathbf{y}_r) =(\mathbf{x}_1,\mathbf{y}_1)\cdots(\mathbf{x}_r,\mathbf{%
y}_r)
\end{equation*}
can be extended by linearity to the inner product on $V^{\otimes r}$ and
then restricted to $\mathsf{V}$.

At microscopic length scales, \emph{spatial randomness} of the material
needs to be taken into account. Mathematically, there is a probability space
$(\Omega,\mathfrak{F},\mathsf{P})$ and a function $\mathsf{T}(A,\omega)\colon%
\mathcal{B}\times\Omega\to\mathsf{V}$ such that for any fixed $A_0\in\mathsf{%
V}$ and for any Borel set $B\subseteq\mathsf{V}$ the inverse image $\mathsf{T%
}^{-1}(A_0,B)$ is an event. The map $\mathsf{T}(\mathbf{x},\omega)$ is a
\emph{random field}.
\end{example}

Translate the whole body $\mathcal{B}$ by a vector $\mathbf{x}\in V$. The
random fields $\mathsf{T}(A+\mathbf{x})$ and $\mathsf{T}(A)$ have the same
finite-dimensional distributions. It is therefore convenient to assume that
there is a random field defined \emph{on all of} $E$ such that its
restriction to $\mathcal{B}$ is equal to $\mathsf{T}(A)$. For brevity,
denote the new field by the same symbol $\mathsf{T}(A)$ (but this time $A\in
E$). The random field $\mathsf{T}(A)$ is \emph{strictly homogeneous}, that
is, the random fields $\mathsf{T}(A+\mathbf{x})$ and $\mathsf{T}(A)$ have
the same finite-dimensional distributions. In other words, for each positive
integer $n$, for each $\mathbf{x}\in V$, and for all distinct places $A_1$,
\dots, $A_n\in E$ the random elements $\mathsf{T}(A_1)\oplus\cdots\oplus%
\mathsf{T}(A_n)$ and $\mathsf{T}(A_1+\mathbf{x})\oplus\cdots\oplus\mathsf{T}%
(A_n+\mathbf{x})$ of the direct sum on $n$ copies of the space $\mathsf{V}$
have the same probability distribution.

Let $K$ be the material symmetry group of the material body $\mathcal{B}$
acting in $V$. The group $K$ is a subgroup of the orthogonal group $\mathrm{O%
}(V)$. For simplicity, we assume that the material is fully symmetric, that
is, $K=\mathrm{O}(V)$. Fix a place $O\in\mathcal{B}$ and identify $E$ with $%
V $ by the map $f$ that maps $A\in E$ to $A-O\in V$. Then $K$ acts in $E$
and rotates the body $\mathcal{B}$ by
\begin{equation*}
g\cdot A=f^{-1}gfA,\qquad g\in\mathrm{O}(V),\quad A\in\mathcal{B}.
\end{equation*}
Let $A_0\in\mathcal{B}$. Under the above action of $K$ the point $A_0$
becomes $g\cdot A_0$. The random tensor $\mathsf{T}(A_0)$ becomes $U(g))%
\mathsf{T}(A_0)$, where $U$ is the restriction of the orthogonal
representation $g\mapsto g^{\otimes r}$ of the group $\mathrm{O}(V)$ to the
subspace $\mathsf{V}$ of the space $V^{\otimes r}$. The random fields $%
\mathsf{T}(g\cdot A)$ and $U(g))\mathsf{T}(A)$ must have the same
finite-dimensional distributions, because $g\cdot A_0$ is the same material
point in a different place. Note that this property does not depend on a
particular choice of the place $O$, because the field is strictly
homogeneous. We call such a field \emph{strictly isotropic}.

Assume that the random field $\mathsf{T}(A)$ is \emph{second-order}, that is
\begin{equation*}
\mathsf{E}[\|\mathsf{T}(A)\|^2]<\infty,\qquad A\in E.
\end{equation*}
Define the \emph{one-point correlation tensor} of the field $\mathsf{T}(A)$
by
\begin{equation*}
\langle\mathsf{T}(A)\rangle=\mathsf{E}[\mathsf{T}(A)]
\end{equation*}
and its \emph{two-point correlation tensor} by
\begin{equation*}
\langle\mathsf{T}(A),\mathsf{T}(B)\rangle=\mathsf{E}[(\mathsf{T}(A) -\langle%
\mathsf{T}(A)\rangle)\otimes(\mathsf{T}(B) -\langle\mathsf{T}(B)\rangle)].
\end{equation*}
Assume that the field $\mathsf{T}(A)$ is \emph{mean-square continuous}, that
is, its two-point correlation tensor $\langle\mathsf{T}(A),\mathsf{T}%
(B)\rangle\colon E\times E\to\mathsf{V}\otimes\mathsf{V}$ is a continuous
function.

Note that \citet{MR3064996} had shown that any finite-variance isotropic
random field on a compact group is necessarily mean-square continuous under
standard measurability assumptions, and hence its covariance function is
continuous. In the related settings, the characterisation of covariance
function for a real homogeneous isotropic random field in $d$-dimensional
Euclidean space was given in the classical paper by \citet{MR1503439}, where
it was conjectured that the only form of discontinuity which could be
allowed for such a function would occur at the origin. This conjecture was
proved by \citet{MR0083534} for $d\geq 2$. This result was widely used in
Geostatistics (see, i.e., \citet{MR1671159}, among the others), who argued
that the homogenous and isotropic random field could be expressed as a
mean-square continuous component and what they called ``nugget effect'',
e.g. a purely discontinuous component. In fact this latter component should
be necessarily non-measurable (see, i.e., \citet[Example~1.2.5]{MR583435}.
The relation between measurability and mean-square continuity in non-compact
situation is still unclear even for scalar random fields. That is why we
assume in this paper that our random fields are mean-square continuous, and
hence their covariance functions are continuous.

If the field $\mathsf{T}(A)$ is strictly homogeneous, then its one-point
correlation tensor is a constant tensor in $\mathsf{V}$, while its two-point
correlation tensor is a function of the vector $B-A$, i.e., a function on $V$%
. Call such a field \emph{wide-sense homogeneous}.

Similarly, if the field $\mathsf{T}(A)$ is strictly isotropic, then we have
\begin{equation}  \label{eq:3}
\begin{aligned} \langle\mathsf{T}(g\cdot
A)\rangle&=U(g)\langle\mathsf{T}(A)\rangle,\\ \langle\mathsf{T}(g\cdot
A),\mathsf{T}(g\cdot B)\rangle &=(U\otimes
U)(g)\langle\mathsf{T}(A),\mathsf{T}(B)\rangle. \end{aligned}
\end{equation}

\begin{definition}
\label{def:1} A random field $\mathsf{T}(A)$ is called \emph{wide-sense
isotropic} if its one-point and two-point correlation tensors satisfy %
\eqref{eq:3}.
\end{definition}

For simplicity, identify the field $\{\,\mathsf{T}(A)\colon A\in E\,\}$
defined on $E$ with the field $\{\,\mathsf{T}^{\prime}(\mathbf{x})\colon%
\mathbf{x}\in V\,\}$ defined by $\mathsf{T}^{\prime}(\mathbf{x})=\mathsf{T}%
(O+\mathbf{x})$. Introduce the Cartesian coordinate system $(x,y,z)$ in $V$.
Use the introduced system to identify $V$ with the coordinate space $\mathbb{%
R}^3$ and $\mathrm{O}(V)$ with $\mathrm{O}(3)$. Call $\mathbb{R}^3$ the
\emph{space domain}. The action of $\mathrm{O}(3)$ on $\mathbb{R}^3$ is the
matrix-vector multiplication.

Definition~\ref{def:1} was used by many authors including \citet{MR0094844}, %
\citet{MR2406668}, \citet{sobczyk2012stochastic}.

There is another definition of isotropy.

\begin{definition}[\citep{MR0094844}]
\label{def:2} A random field $\mathsf{T}(A)$ is called a \emph{%
multidimensional scalar wide-sense isotropic} if its one-point correlation
tensor is a constant, while the two-point correlation tensor $\langle\mathsf{%
T}(\mathbf{x},\mathsf{T}(\mathbf{y})\rangle$ depends only on $\|\mathbf{y}-%
\mathbf{x}\|$.
\end{definition}

It is easy to see that Definition~\ref{def:2} is a particular case of
Definition~\ref{def:1} when the representation $U$ is trivial, that is, maps
all elements $g\in K$ to the identity operator.

In the case of $r=0$, the complete description of the two-point correlation
functions of scalar homogeneous and isotropic random fields is as follows.
Recall that a measure $\mu$ defined on the Borel $\sigma$-field of a
Hausdorff topological space $X$ is called \emph{Borel measure}.

\begin{theorem}
\label{th:1} Formula
\begin{equation}  \label{eq:1}
\langle T(\mathbf{x}),T(\mathbf{y})\rangle=\int^{\infty}_0\frac{%
\sin(\lambda\|\mathbf{y}-\mathbf{x}\|)} {\lambda\|\mathbf{y}-\mathbf{x}\|}\,%
\mathrm{d}\mu(\lambda)
\end{equation}
establishes a one-to-one correspondence between the set of two-point
correlation functions of homogeneous and isotropic random fields $T(\mathbf{x%
})$ on the space domain $\mathbb{R}^3$ and the set of all finite Borel
measures $\mu$ on the interval $[0,\infty)$.
\end{theorem}

Theorem~\ref{th:1} is a translation of the result proved by \citet{MR1503439}
to the language of random fields. This translation is performed as follows.
Assume that $B(\mathbf{x})$ is a two-point correlation function of a
homogeneous and isotropic random field $T(\mathbf{x})$. Let $n$ be a
positive integer, let $\mathbf{x}_1$, \dots, $\mathbf{x}_n$ be $n$ distinct
points in $\mathbb{R}^3$, and let $c_1$, \dots, $c_n$ be $n$ complex
numbers. Consider the random variable $X=\sum^n_{j=1}c_j[T(\mathbf{x}%
_j)-\langle T(\mathbf{x}_j)\rangle]$. Its variance is non-negative:
\begin{equation*}
\mathsf{E}[X^2]=\sum_{j,k=1}^{n}c_j\overline{c_k}\langle T(\mathbf{x}_j),T(%
\mathbf{x}_k)\rangle\geq 0.
\end{equation*}
In other words, the two-point correlation function $\langle T(\mathbf{x}),T(%
\mathbf{y})\rangle$ is a non\-ne\-ga\-tive-de\-fi\-nite function. Moreover, it is
continuous, because the random field $T(\mathbf{x})$ is mean-square
continuous, and depends only on the distance $\|\mathbf{y}-\mathbf{x}\|$
between the points $\mathbf{x}$ and $\mathbf{y}$, because the field is
homogeneous and isotropic. \citet{MR1503439} proved that Equation~%
\eqref{eq:1} describes all of such functions.

Conversely, assume that the function $\langle T(\mathbf{x}),T(\mathbf{y}%
)\rangle$ is described by Equation~\eqref{eq:1}. The centred Gaussian random
field with the two-point correlation function \eqref{eq:1} is homogeneous
and isotropic. In other words, there is a link between the theory of random
fields and the theory of positive-definite functions.

In what follows, we consider the fields with absolutely continuous spectrum.

\begin{definition}[\citet{MR1009786}]
A homogeneous and isotropic random field $T(\mathbf{x})$ has an \emph{%
absolutely continuous spectrum} if the measure $\mu$ is absolutely
continuous with respect to the measure $4\pi\lambda^2\,\mathrm{d}\lambda$,
i.e., there exist a nonnegative measurable function $f(\lambda)$ such that
\begin{equation*}
\int^{\infty}_0\lambda^2f(\lambda)\,\mathrm{d}\lambda<\infty
\end{equation*}
and $d\mu(\lambda)=4\pi\lambda^2f(\lambda)\,\mathrm{d}\lambda$. The function
$f(\lambda)$ is called the \emph{isotropic spectral density} of the random
field $T(\mathbf{x})$.
\end{definition}

\begin{example}[The Mat\'{e}rn two-point correlation function]
\label{ex:2} Consider a two-point correlation function of a scalar random
field $T(\mathbf{x})$ of the form
\begin{equation}  \label{eq:2}
\left\langle T(\mathbf{x}),T(\mathbf{y})\right\rangle =M_{\nu ,a}\left(
\mathbf{x},\mathbf{y}\right) =\frac{2^{1-\nu }\sigma ^{2}}{\Gamma \left( \nu
\right) }\left( a\left\Vert \mathbf{x}-\mathbf{y}\right\Vert \right) ^{{}\nu
}K_{{}\nu }\left( a\left\Vert \mathbf{x}-\mathbf{y}\right\Vert \right) ,\quad
\end{equation}%
where $\sigma ^{2}>0,a>0,\nu >0$ and $K_{{}\nu }\left( z\right) $ is the
Bessel function of the third kind of order $\nu$. Here, the parameter $\nu$
measures the differentiability of the random field; the parameter $\sigma $
is its variance and the parameter $a$ measures how quickly the correlation
function of the random field decays with distance. The corresponding
isotropic spectral density is
\begin{equation*}
f\left(\lambda\right) =f_{\nu ,a,\sigma ^{2}}\left(\lambda\right) =\frac{%
\sigma ^{2}\Gamma \left( \nu +\frac{3}{2}\right) a^{2\nu }}{2\pi
^{3/2}\left( a^{2}+\lambda^{2}\right) ^{\nu +\frac{3}{2}}},\quad \lambda\geq
0.
\end{equation*}
\end{example}

Note that Example~\ref{ex:2} demonstrates another link, this time between
the theory of random fields and the theory of special functions.

In this paper, we consider the following problem. How to define the Mat\'{e}%
rn two-point correlation tensor for the case of $r>0$? A particular answer
to this question can be formulated as follows.

\begin{example}[Parsimonious Mat\'{e}rn model, \citet{MR2752612}]
\label{ex:3} We assume that the vector random field
\begin{equation*}
T\left( \mathbf{x}\right) =\left( T_{1}\left( \mathbf{x}\right)
,\dots,T_{m}\left(\mathbf{x}\right) \right)^{\top} ,\qquad\mathbf{x}\in
\mathbb{R}^{3}
\end{equation*}
has the two-point correlation tensor $B\left(\mathbf{x},\mathbf{y}\right)
=\left( B_{ij}\left( \mathbf{x},\mathbf{y}\right) \right) _{1\leq i,j\leq
m}. $ It is not straightforward to specify the cross-covariance functions $%
B_{ij}\left( \mathbf{x}\right) ,1\leq i,j\leq m,i\neq j$ as non-trivial,
valid parametric models\ because of the requirement of their non-negative
definiteness. In the multivariate Mat\'{e}rn model, each marginal covariance
function%
\begin{equation*}
B_{ii}\left( \mathbf{x},\mathbf{y}\right) =\sigma _{i}^{2}M_{\nu
_{i},a_{i}}\left(\mathbf{x},\mathbf{y}\right) ,i=1,...,m,
\end{equation*}%
is of the type \eqref{eq:2} with the isotropic spectral density $%
f_{ii}(\lambda)=f_{\nu _{i},a_{i},\sigma _{i}^{2}}\left( \lambda\right) .$

Each cross-covariance function%
\begin{equation*}
B_{ij}\left(\mathbf{x},\mathbf{y}\right) =B_{ji}\left(\mathbf{x},\mathbf{y}%
\right) =b_{ij}\sigma _{i}\sigma _{j}M_{\nu _{ij},a_{ij}}\left( \mathbf{x},%
\mathbf{y}\right) ,1\leq i,j\leq m,i\neq j
\end{equation*}%
is also a Mat\'{e}rn function with co-location correlation coefficient $%
b_{ij},$ smoothness parameter $\nu _{ij}$ and scale parameter $a_{ij}.$The
spectral densities are
\begin{equation*}
f_{ij}\left( \mathbf{x}\right) =f_{\nu _{ij},a_{ij},,b_{ij}\sigma _{i}\sigma
_{j}}\left( \mathbf{x}\right) ,1\leq i,j\leq m,i\neq j.
\end{equation*}

The question then is to determine the values of $\nu _{ij},a_{ij}$ and $%
b_{ij}$ so that the non-negative definiteness condition is satisfied. Let $%
m\geq 2$. Suppose that%
\begin{equation*}
\nu _{ij}=\frac{1}{2}\left( \nu _{i}+\nu _{j}\right) ,1\leq i,j\leq m,i\neq
j,
\end{equation*}%
and that there is a common scale parameter in the sense that there exists an
$a>0$ such that
\begin{equation*}
a_{i}=...=a_{m}=a,\text{ and }a_{ij}=a\text{ for }1\leq i,j\leq m,i\neq j.
\end{equation*}%
Then the multivariate Mat\'{e}rn model provides a valid second-order
structure in $\mathbb{R}^{3}$ if%
\begin{equation*}
b_{ij}=\beta _{ij}\left[ \frac{\Gamma \left( \nu _{i}+\frac{3}{2}\right) }{%
\Gamma \left( \nu _{i}\right) }\frac{\Gamma \left( \nu _{j}+\frac{3}{2}%
\right) }{\Gamma \left( \nu _{j}\right) }\right] ^{1/2}\frac{\Gamma \left(
\frac{1}{2}\left( \nu _{i}+\nu _{j}\right) \right) }{\Gamma \left( \frac{1}{2%
}\left( \nu _{i}+\nu _{j}\right) +\frac{3}{2}\right) }
\end{equation*}%
for $1\leq i,j\leq m,i\neq j,$ where the matrix $\left( \beta _{ij}\right)
_{i,j=1,...,m}$ has diagonal elements $\beta _{ii}=1$ for $i=1,...,m,$ and
off-diagonal elements $\beta _{ij},1\leq i,j\leq m,i\neq j$ so that it is
symmetric and non-negative definite. \newline
\end{example}

\begin{example}[Flexible Mat\'{e}rn model]
Consider the vector random field $\mathsf{T}(\mathbf{x})\in \mathbb{R}^{m},%
\mathbf{x}\in \mathbb{R}^{3}$ with the two-point covariance tensor%
\begin{equation*}
\left\langle T_{i}(x),T_{j}(\mathbf{y})\right\rangle =B_{ij}(\mathbf{x},%
\mathbf{y})=\bar{B}_{ij}(\mathbf{y}-\mathbf{x})=\sigma _{ij}M_{\nu
_{ij},a_{ij}}\left( \mathbf{x},\mathbf{y}\right) ,1\leq i,j\leq m,
\end{equation*}%
where again
\begin{equation*}
M_{\nu ,a}\left( \mathbf{x},\mathbf{y}\right) =\frac{2^{1-\nu }\sigma ^{2}}{%
\Gamma \left( \nu \right) }\left( a\left\Vert \mathbf{y}-\mathbf{x}%
\right\Vert \right) ^{{}\nu }K_{{}\nu }\left( a\left\Vert \mathbf{y}-\mathbf{%
x}\right\Vert \right) .\quad
\end{equation*}%
We assume that the matrix $\Sigma =(\sigma _{ij})_{1\leq i,j\leq m}=(\sigma
_{ij})>0$ (nonnegative definite), and we denote $\sigma _{i}^{2}=\sigma
_{ii} $, $i=1$, \dots , $m$.

Then the spectral density $F=(f_{ij})_{1\leq i,j\leq m}$ has the entries%
\begin{equation*}
\begin{aligned} f_{ij}(\bm{\lambda} )&=\frac{1}{(2\pi
)^{3}}\int_{\mathbb{R}^{3}}e^{-\mathrm{i}(\bm{\lambda}
,\mathbf{h})}\bar{B}_{ij}(\mathbf{h})\,\mathrm{d}\mathbf{h}\\ &=\sigma
_{ij}a_{ij}^{2\nu _{ij}}\frac{1}{(a_{ij}+\left\Vert \bm{\lambda} \right\Vert
^{2})^{\nu _{ij}+\frac{3}{2}}}\frac{\Gamma (\nu _{ij}+\frac{3}{2})}{\Gamma
(\nu _{ij})},1\leq i,j\leq m,\lambda \in \mathbb{R}^{3}. \end{aligned}
\end{equation*}

We need to find some conditions on parameters $a_{ij}>0,\nu _{ij}>0,$ under
which the matrix $F>0$ (nonnegative definite). The general conditions can be
found in \citet{MR2946043} and \citet{MR2949350}.

Recall that a symmetric, real $m\times \ m$ matrix $\Theta =(\theta
_{ij})_{1\leq i,j\leq m},$ is said to be conditionally negative definite %
\citep{MR1449393}, if the inequality%
\begin{equation*}
\sum_{i=1}^{m}\sum_{j=1}^{m}c_{i}c_{j}\theta _{ij}\leq 0
\end{equation*}
holds for any real numbers $c_{1},...,$ $c_{m,}$ subject to
\begin{equation*}
\sum_{i=1}^{m}c_{i}=0.
\end{equation*}

In general, a necessary condition for the above inequality is%
\begin{equation*}
\theta_{ii}+\theta_{jj}\leq 2\theta_{ij},\qquad i,j=1,...,m,
\end{equation*}
which implies that all entries of a conditionally negative definite matrix
are nonnegative whenever its diagonal entries are non-negative. If all its
diagonal entries vanish, a conditionally negative definite matrix is also
named a Euclidean distance matrix. It is known that $\Theta =(\theta
_{ij})_{1\leq i,j\leq m}$ is conditionally negative definite if and only if
an $m\times \ m$ matrix $S$ with entries $\exp \{-\theta_{ij}u\}$ is
positive definite, for every fixed $u\geq 0$ (cf.
\citet[Theorem
4.1.3]{MR1449393}), or $S=e^{-u\Theta },$ where $e^{\Lambda }$ is an Hadamar
exponential of a matrix $\Lambda .$

Some simple examples of conditionally negative definite matrices are

(i) $\theta_{ij}=\theta_{i}+\theta_{j};$

(ii) $\theta_{ij}=\mathrm{const};$

(iii) $\theta_{ij}=\left\vert\theta_{i}-\theta_{j}\right\vert ;$

(iv) $\theta_{ij}=\left\vert\theta_{i}-\theta_{j}\right\vert ^{2}$

(v) $\theta_{ij}=\max \{\theta_{i},\theta_{j}\};$

(vi) $\theta_{ij}=-\theta_{i}\theta_{j}.$

Recall that the Hadamard product of two matrices $A$ and $B$ is the matrix $%
A\circ B=(A_{ij}\cdot B_{ij})_{1\leq i,j\leq m}.$ By Schur theorem if $%
A>0,B>0,$ then so is $A\circ B.$

Then
\begin{equation*}
F=\Sigma \circ A\circ B\circ C,
\end{equation*}%
where one need to find conditions under which
\begin{equation*}
\begin{aligned} A&=\left( \frac{1}{(1+\left\Vert \bm{\lambda} \right\Vert
^{2}/a_{ij}^{2})^{\nu _{ij}+\frac{3}{2}}}\right) _{1\leq i,j\leq m}\geq
0,\qquad B=\left( \frac{1}{a_{ij}^{3}}\right) _{1\leq i,j\leq m}\geq 0,\\
C&=\left( \frac{\Gamma (\nu _{ij}+\frac{3}{2})}{\Gamma (\nu _{ij})}\right)
_{1\leq i,j\leq m}\geq 0. \end{aligned}
\end{equation*}

We consider first the case 1, in which we assume that
\begin{equation*}
a_{i}=...=a_{m}=a,\text{ }1\leq i,j\leq m.
\end{equation*}%
Then
\begin{equation*}
A=e^{-\frac{3}{2}}\left( \exp \{-\nu _{ij}\log (1+\frac{\left\Vert %
\bm{\lambda}\right\Vert ^{2}}{a^{2}})\}\right) _{1\leq i,j\leq m}\geq 0,
\end{equation*}%
if and only if \ the matrix
\begin{equation*}
Y=\left( -\nu _{ij}\right) _{1\leq i,j\leq m}
\end{equation*}%
is conditionally negative definite (see above examples (i)-(vi)), then for
such $\left( -\nu _{ij}\right) _{1\leq i,j\leq m},$ we have to check that
the matrix $C=(\Gamma (\nu _{ij}+\frac{3}{2})/\Gamma (\nu _{ij})_{1\leq
i,j\leq m}\geq 0.$ This class is not empty, since it is included the case of
the so-called parsimonious model: $\nu _{ij}=\frac{\nu _{i}+\nu _{j}}{2}$
(see Example~\ref{ex:3}). Thus, for the case 1, the following multivariate
Mat\'{e}rn models are valid under the following conditions (see, %
\citet{MR2946043,MR2949350}):

A1) Assume that

i) $a_{i}=...=a_{m}=a,$ $1\leq i,j\leq m;$

ii) $-\nu _{ij}$ ,$1\leq i,j\leq m;$ form a conditionally non-negative
matrices;

iii) $\sigma _{ij}\frac{\Gamma (\nu _{ij}+\frac{3}{2})}{\Gamma (\nu _{ij})}%
,1\leq i,j\leq m;$ form a non-negative definite matrices.

Consider the case 2:
\begin{equation*}
\nu _{ij}=\nu >0,\text{ }1\leq i,j\leq m.
\end{equation*}
Then the following multivariate Mat\'{e}rn models are valid under the
following conditions \citep{MR2949350}:

A2) either

a) $-a_{ij}^{2}$ ,$1\leq i,j\leq m,$ form a conditionally non-negative
matrix and $\sigma _{ij}a_{ij}^{2\nu },1\leq i,j\leq m,$ form non-negative
definite matrices;

or

b) $-a_{ij}^{-2}$ ,$1\leq i,j\leq m,$ form a conditionally non-negative
matrix and $\sigma _{ij}/a_{ij}^{3},1\leq i,j\leq m,$ form non-negative
definite matrices.

These classes of Mat\'{e}rn models are not empty since in the case of
parsimonious model they are consistent with \citet[Theorem~1]{MR2752612}.
For the parsimonious model form this paper $($ $\nu _{ij}=\frac{\nu
_{ii}+\nu _{jj}}{2},1\leq i,j\leq m),$ the following multivariate Mat\'{e}rn
models are valid under conditions

A3) either

a) $\ \nu _{ij}=\frac{\nu _{ii}+\nu _{jj}}{2},a_{ij}^{2}=\frac{%
a_{ii}^{2}+a_{jj}^{2}}{2},1\leq i,j\leq m,$ and $\sigma _{ij}a_{ij}^{2\nu
_{ij}}/\Gamma (\nu _{ij}),1\leq i,j\leq m,$form non-negative definite
matrices;

or

b) $\nu _{ij}=\frac{\nu _{ii}+\nu _{jj}}{2},a_{ij}^{-2}=\frac{%
a_{ii}^{-2}+a_{jj}^{-2}}{2},1\leq i,j\leq m,$and $\sigma
_{ij}/a_{ij}^{3}/\Gamma (\nu _{ij}),1\leq i,j\leq m,$ form non-negative
definite matrices;

The most general conditions and new examples can be found in %
\citet{MR2946043} and \citet{MR2949350}. The paper by \citet{MR3353096}
reviews the main approaches to building multivariate correlation and
covariance structures, including the multivariate Mat\'{e}rn models.
\end{example}

\begin{example}[Dual Mat\'{e}rn models]
\label{ex:dual} Adapting the so-called duality theorem (see, i.e., %
\citet{MR2557625}), one can show that under the conditions A1, A2 or A3
\begin{equation*}
\frac{1}{(1+\left\Vert \mathbf{h}\right\Vert ^{2})^{\nu _{ij}+\frac{3}{2}}}%
=\int_{\mathbb{R}^{3}}e^{\mathrm{i}(\bm{\lambda },\mathbf{h})}s_{ij}(%
\bm{\lambda })d\bm{\lambda },\qquad 1\leq i,j\leq m,
\end{equation*}%
where%
\begin{equation*}
s_{ij}(\bm{\lambda )=}\frac{1}{(2\pi )^{3}2^{\nu _{ij}-1}\Gamma (\nu _{ij}+%
\frac{3}{2})}(\left\Vert \bm{\lambda }\right\Vert )^{\nu _{ij}}K_{\nu
_{ij}}(\left\Vert \bm{\lambda }\right\Vert ),\qquad\bm{\lambda \in }\mathbb{R%
}^{3},1\leq i,j\leq m,
\end{equation*}%
is the valid spectral density of the vector random field with correlation
structure $((1+\left\Vert \mathbf{h}\right\Vert ^{2})^{-(\nu _{ij}+\frac{3}{2%
})})_{1\leq i,j\leq m}=(D_{ij}(\mathbf{h}))_{1\leq i,j\leq m}$. We will call
it the \emph{dual Mat\'{e}rn model}.

Note that for the Mat\'{e}rn models
\begin{equation*}
\int_{\mathbb{R}^{3}}\bar{B}_{ij}(\mathbf{x})d\mathbf{x}<\infty .
\end{equation*}

\bigskip This condition is known as short range dependence, while for the
dual Mat\'{e}rn model, the long range dependence is possible:

\begin{equation*}
\int_{\mathbb{R}^{3}}D_{ij}(\mathbf{h})d\mathbf{h}=\infty ,\text{ if }0<\nu
_{ij}<\frac{3}{2}.
\end{equation*}
\end{example}

When $m=3$, the random field of Example~\ref{ex:3} is scalar isotropic but
not isotropic. How to construct examples of homogeneous and \emph{isotropic}
vector and tensor random fields with Mat\'{e}rn two-point correlation
tensors?

To solve this problem, we develop a sketch of a general theory of
homogeneous and isotropic tensor-valued random fields in Section~\ref%
{sec:general}. This theory was developed by \citet{MR3336288,MR3493458}. In
particular, we explain another two links: one leads from the theory of
random fields to classical invariant theory, another one was established
recently and leads from the theory of random fields to the theory of convex
compacta.

In Section~\ref{sec:examples}, we give examples of Mat\'{e}rn homogeneous
and isotropic tensor-valued random fields. Finally, in Appendices we shortly
describe the mathematical terminology which is not always familiar to
specialists in probability: tensors, group representations, and classical
invariant theory. For different aspects of theory of random fields see also %
\citet{MR1687092} and \citet{MR2870527}.

\section{A sketch of a general theory}

\label{sec:general}

Let $r$ be a nonnegative integer, let $\mathsf{V}$ be an invariant subspace
of the representation $g\mapsto g^{\otimes r}$ of the group $\mathrm{O}(3)$,
and let $U$ be the restriction of the above representation to $\mathsf{V}$.
Consider a homogeneous $\mathsf{V}$-valued random field $\mathsf{T}(\mathbf{x%
})$, $\mathbf{x}\in\mathbb{R}^3$. Assume it is isotropic, that is, satisfies %
\eqref{eq:3}. It is very easy to see that its one-point correlation tensor $%
\langle\mathsf{T}(\mathbf{x})\rangle$ is an arbitrary element of the
isotypic subspace of the space $\mathsf{V}$ that corresponds to the trivial
representation. In particular, in the case of $r=0$ the representation $U$
is trivial, and $\langle\mathsf{T}(\mathbf{x})\rangle$ is an arbitrary real
number. In the case of $r=1$ we have $U(g)=g$. This representation does not
contain a trivial component, therefore $\langle\mathsf{T}(\mathbf{x})\rangle=%
\mathbf{0}$. In the case of $r=2$ and $U(g)=\mathsf{S}^2(g)$ the isotypic
subspace that corresponds to the trivial representation is described in
Example~\ref{ex:8}, we have $\langle\mathsf{T}(\mathbf{x})\rangle=CI$, where
$C$ is an arbitrary real number, and $I$ is the identity operator in $%
\mathbb{R}^3$, and so on.

Can we quickly describe the two-point correlation tensor in the same way?
The answer is positive. Indeed, the second equation in \eqref{eq:3} means
that $\langle\mathsf{T}(\mathbf{x}),\mathsf{T}(\mathbf{y})\rangle$ is a
measurable covariant of the pair $(g,U)$. The integrity basis for polynomial
invariants of the defining representation contains one element $I_1=\|%
\mathbf{x}\|^2$. By the Wineman--Pipkin theorem (Appendix~\ref{ap:tensors},
Theorem~\ref{th:Wineman-Pipkin}), we obtain
\begin{equation*}
\langle\mathsf{T}(\mathbf{x}),\mathsf{T}(\mathbf{y})\rangle=\sum_{l=1}^{L}
\varphi_l(\|\mathbf{y}-\mathbf{x}\|^2)\mathsf{T}_l(\mathbf{y}-\mathbf{x}),
\end{equation*}
where $\mathsf{T}_l(\mathbf{y}-\mathbf{x})$ are the basic covariant tensors
of the representation $U$.

For example, when $r=1$, the basis covariant tensors of the defining
representations are $\delta_{ij}$ and $x_ix_j$ by the result of %
\citet{MR1488158} mentioned in Appendix~\ref{ap:invariant}. We obtain the
result by \citet{MR0001702}:
\begin{equation*}
\langle\mathsf{T}(\mathbf{x}),\mathsf{T}(\mathbf{y})\rangle= \varphi_1(\|%
\mathbf{y}-\mathbf{x}\|^2)\delta_{ij} +\varphi_2(\|\mathbf{y}-\mathbf{x}\|^2)%
\frac{(y_i-x_i)(y_j-x_j)}{\|\mathbf{y}-\mathbf{x}\|^2}.
\end{equation*}

When $r=2$ and $U(g)=\mathsf{S}^2(g)$, the three rank~$4$ isotropic tensors
are $\delta_{ij}\delta_{kl}$, $\delta_{ik}\delta_{jl}$, and $%
\delta_{il}\delta_{jk}$. Consider the group $\Sigma$ of order~$8$ of the
permutations of symbols $i$, $j$, $k$, and $l$, generated by the
transpositions $(ij)$, $(kl)$, and the product $(ik)(jl)$. The group $\Sigma$
acts on the set of rank~$4$ isotropic tensors and has two orbits. The sums
of elements on each orbit are basis isotropic tensors:
\begin{equation*}
L^1_{ijkl}=\delta_{ij}\delta_{kl},\qquad L^2_{ijkl}=\delta_{ik}\delta_{jl}
+\delta_{il}\delta_{jk}.
\end{equation*}
Consider the case of degree~$2$ and of order~$4$. For the pair of representations $%
(g^{\otimes 4},(\mathbb{R}^3)^{\otimes 4})$ and $(g,\mathbb{R}^3)$ we have $%
6 $~covariant tensors:
\begin{equation*}
\delta_{il}x_jx_k,\delta_{jk}x_ix_{l},\delta_{jl}x_ix_k,
\delta_{ik}x_jx_{l},\delta_{kl}x_ix_j,\delta_{ij}x_kx_{l}.
\end{equation*}
The action of the group $\Sigma$ has $2$~orbits, and the symmetric covariant
tensors are
\begin{equation*}
\begin{aligned} \|\mathbf{x}\|^2L^3_{ijkl}(\mathbf{x})&=\delta_{il}x_jx_k
+\delta_{jk}x_ix_{l}+\delta_{jl}x_ix_k+\delta_{ik}x_jx_{l},\\
\|\mathbf{x}\|^2L^4_{ijkl}(\mathbf{x})&=\delta_{kl}x_ix_j
+\delta_{ij}x_kx_{l}. \end{aligned}
\end{equation*}
In the case of degree~$4$ and of order~$4$ we have only one covariant:
\begin{equation*}
\|\mathbf{x}\|^4L^5_{ijkl}(\mathbf{x})=x_ix_jx_kx_{l}.
\end{equation*}
The result by \citet{Lomakin1964}
\begin{equation*}
\langle\mathsf{T}(\mathbf{x}),\mathsf{T}(\mathbf{y})\rangle
=\sum_{m=1}^{5}\varphi_m(\|\mathbf{y}-\mathbf{x}\|^2) L^m_{ijkl}(\mathbf{y}-%
\mathbf{x})
\end{equation*}
easily follows.

The case of $r=3$ will be considered in details elsewhere.

When $r=4$ and $U(g)=\mathsf{S}^2(S^2(g))$, the situation is more delicate.
There are $8$ symmetric isotropic tensors connected by $1$ syzygy, $13$
basic covariant tensors of degree~$2$ and of order~$8$ connected by $3$%
~syzygies, $10$ basic covariant tensors of degree~$4$ and of order~$8$
connected by $2$~syzygies, $3$ basic covariant tensors of degree~$6$ and of
order~$8$, and $1$ basic covariant tensor of degree~$8$ and of order~$8$,
see \citet{Malyarenko2016a} and \citet{Malyarenko2016} for details. It follows that there are $29$
independent basic covariant tensors. The result by \citet{Lomakin1965}
includes only $15$ of them and is therefore incomplete.

How to find the functions $\varphi_m$? In the case of $r=0$, the answer is
given by Theorem~\ref{th:1}:
\begin{equation*}
\varphi_1(\|\mathbf{y}-\mathbf{x}\|^2)=\int^{\infty}_0 \frac{\sin(\lambda\|%
\mathbf{y}-\mathbf{x}\|)} {\lambda\|\mathbf{y}-\mathbf{x}\|}\,\mathrm{d}%
\mu(\lambda).
\end{equation*}
In the case of $r=1$, the answer has been found by \citet{MR0094844}:
\begin{equation}  \label{eq:Yaglom}
\begin{aligned}
\varphi_1(\|\mathbf{y}-\mathbf{x}\|^2)&=\frac{1}{\rho^2}\left(\int^{%
\infty}_0j_2(\lambda\rho)
\,\mathrm{d}\Phi_2(\lambda)-\int^{\infty}_0j_1(\lambda\rho)
\,\mathrm{d}\Phi_1(\lambda)\right),\\
\varphi_2(\|\mathbf{y}-\mathbf{x}\|^2)&=\int^{\infty}_0\frac{j_1(\lambda%
\rho)}{\lambda\rho} \,\mathrm{d}\Phi_1(\lambda)
+\int^{\infty}_0\left(j_0(\lambda\rho)-\frac{j_1(\lambda\rho)}{\lambda\rho}
\right)\,\mathrm{d}\Phi_2(\lambda), \end{aligned}
\end{equation}
where $\rho=\|\mathbf{y}-\mathbf{x}\|$, $j_n$ are the spherical Bessel
functions, and $\Phi_1$ and $\Phi_2$ are two finite measures on $[0,\infty)$
with $\Phi_1(\{0\})=\Phi_2(\{0\})$.

In the general case, we proceed in steps. The main idea is simple. We
describe all homogeneous random fields and throw away those that are not
isotropic. The homogeneous random fields are described by the following
result.

\begin{theorem}
\label{th:Kolmogorov} Formula
\begin{equation}  \label{eq:9}
\langle\mathsf{T}(\mathbf{x}),\mathsf{T}(\mathbf{y})\rangle =\int_{\hat{%
\mathbb{R}}^3}e^{\mathrm{i}(\mathbf{p},\mathbf{y}-\mathbf{x})} \,\mathrm{d}%
\mu(\mathbf{p})
\end{equation}
establishes a one-to-one correspondence between the set of the two-point
correlation tensors of homogeneous random fields $\mathsf{T}(\mathbf{x})$ on
the \emph{space domain} $\mathbb{R}^3$ with values in a \emph{complex}
finite-dimensional space $\mathsf{V}_{\mathbb{C}}$ and the set of all
measures $\mu$ on the Borel $\sigma$-field $\mathfrak{B}(\hat{\mathbb{R}}^3)$
of the \emph{wavenumber domain} $\hat{\mathbb{R}}^3$ with values in the cone
of nonnegative-definite Hermitian operators in $\mathsf{V}_{\mathbb{C}}$.
\end{theorem}

This theorem was proved by \citep{MR0003440,MR0003441} for one-dimensional
stochastic processes. Kolmogorov's results have been further developed by %
\citet{MR0006609}, \citet{MR0013259}, \citep{MR0015712,MR0015713} among
others.

We would like to write as many formulae as possible in a coordinate-free
form, like \eqref{eq:9}. To do that, let $J$ be a \emph{real structure} in
the space $\mathsf{V}_{\mathbb{C}}$, that is, a map $j\colon\mathsf{V}_{%
\mathbb{C}}\to\mathsf{V}_{\mathbb{C}}$ with

\begin{itemize}
\item $J(\mathsf{x}+\mathsf{y})=J(\mathsf{x})+J(\mathsf{y})$, $\mathsf{x}$, $%
\mathsf{y}\in\mathsf{V}_{\mathbb{C}}$.

\item $J(\alpha\mathsf{x})=\overline{\alpha}J(\mathsf{x})$, $\mathsf{x}\in%
\mathsf{V}_{\mathbb{C}}$, $\alpha\in\mathbb{C}$.

\item $J(J(\mathsf{x}))=\mathsf{x}$, $\mathsf{x}\in\mathsf{V}_{\mathbb{C}}$.
\end{itemize}

Any tensor $\mathsf{x}\in\mathsf{V}_{\mathbb{C}}$ can be written as $\mathsf{%
x}=\mathsf{x}^++\mathsf{x}^-$, where
\begin{equation*}
\mathsf{x}^+=\frac{1}{2}(\mathsf{x}+J\mathsf{x}),\qquad \mathsf{x}^-=\frac{1%
}{2}(\mathsf{x}-J\mathsf{x}).
\end{equation*}
Denote
\begin{equation*}
\mathsf{V}^+=\{\,\mathsf{x}\in\mathsf{V}_{\mathbb{C}}\colon J\mathsf{x}=%
\mathsf{x}\,\},\qquad\mathsf{V}^-=\{\,\mathsf{x}\in\mathsf{V}_{\mathbb{C}%
}\colon J\mathsf{x}=-\mathsf{x}\,\}.
\end{equation*}
Both sets $\mathsf{V}^+$ and $\mathsf{V}^-$ are real vector spaces. If the
values of the random field $\mathsf{T}(\mathbf{x})$ lie in $\mathsf{V}^+$,
then the measure $\mu$ satisfies the condition
\begin{equation}  \label{eq:10}
\mu(-A)=\mu^{\top}(A)
\end{equation}
for all Borel subsets $A\subseteq\hat{\mathbb{R}}^3$, where $-A=\{\,-\mathbf{%
p}\colon\mathbf{p}\in A\,\}$.

Next, the following Lemma can be proved. Let $\mathbf{p}=(\lambda,\varphi_{%
\mathbf{p}},\theta_{\mathbf{p}})$ be the spherical coordinates in the
wavenumber domain.

\begin{lemma}
A homogeneous random field described by \emph{\eqref{eq:9}} and \emph{%
\eqref{eq:10}} is isotropic if and only if its two-point correlation tensor
has the form
\begin{equation}  \label{eq:12}
\langle\mathsf{T}(\mathbf{x}),\mathsf{T}(\mathbf{y})\rangle=\frac{1}{4\pi}
\int_{0}^{\infty}\int_{S^2}e^{\mathrm{i}(\mathbf{p},\mathbf{y}-\mathbf{x})}
f(\lambda,\varphi_{\mathbf{p}},\theta_{\mathbf{p}})\sin\theta_{\mathbf{p}} \,%
\mathrm{d}\varphi_{\mathbf{p}}\,\mathrm{d}\theta_{\mathbf{p}}\,\mathrm{d}
\nu(\lambda),
\end{equation}
where $\nu$ is a finite measure on the interval $[0,\infty)$, and where $f$
is a measurable function taking values in the set of all symmetric
nonnegative-definite operators on $\mathsf{V}^+$ with unit trace and
satisfying the condition
\begin{equation}  \label{eq:11}
f(g\mathbf{p})=\mathsf{S}^2(U)(g)f(\mathbf{p}),\qquad\mathbf{p}\in\hat{%
\mathbb{R}}^3, \quad g\in\mathrm{O}(3).
\end{equation}
\end{lemma}

When $\lambda=0$, condition~\eqref{eq:11} gives $f(\mathbf{0})=\mathsf{S}%
^2(U)(g)f(\mathbf{0})$ for all $g\in\mathrm{O}(3)$. In other words, the
tensor $f(\mathbf{0})$ lies in the isotypic subspace of the space $\mathsf{S}%
^2(\mathsf{V^+})$ that corresponds to the trivial representation of the
group $\mathrm{O}(3)$, call it $\mathsf{H}_1$. The intersection of $\mathsf{H%
}_1$ with the set of all symmetric nonnegative-definite operators on $%
\mathsf{V}^+$ with unit trace is a convex compact set, call it $\mathcal{C}%
_1 $.

When $\lambda>0$, condition~\eqref{eq:11} gives $f(\lambda,0,0)=\mathsf{S}%
^2(U)(g)f(\lambda,0,0)$ for all $g\in\mathrm{O}(2)$, because $\mathrm{O}(2)$
is the subgroup of $\mathrm{O}(3)$ that fixes the point $(\lambda,0,0)$. In
other words, consider the restriction of the representation $\mathsf{S}^2(U)$
to the subgroup $\mathrm{O}(2)$. The tensor $f(\lambda,0,0)$ lies in the
isotypic subspace of the space $\mathsf{S}^2(\mathsf{V^+})$ that corresponds
to the trivial representation of the group $\mathrm{O}(2)$, call it $\mathsf{%
H}_0$. We have $\mathsf{H}_1\subset\mathsf{H}_0$, because $\mathrm{O}(2)$ is
a subgroup of $\mathrm{O}(3)$. The intersection of $\mathsf{H}_0$ with the
set of all symmetric nonnegative-definite operators on $\mathsf{V}^+$ with
unit trace is a convex compact set, call it $\mathcal{C}_0$.

Fix an orthonormal basis $\mathsf{T}^{0,1,0}$, \dots, $\mathsf{T}^{0,n_0,0}$
of the space $\mathsf{H}_1$. Assume that the space $\mathsf{H}_0\ominus%
\mathsf{H}_1$ has the non-zero intersection with the spaces of $n_1$ copies
of the irreducible representation $U^{2g}$, $n_2$ copies of the irreducible
representation $U^{4g}$, \dots, $n_r$ copies of the irreducible
representation $U^{2rg}$ of the group $\mathrm{O}(3)$, and let $\mathsf{T}%
^{2\ell,n,m}$, $-2\ell\leq m\leq 2\ell$, be the tensors of the Gordienko
basis of the $n$th copy of the representation $U^{2\ell g}$. We have
\begin{equation}  \label{eq:15}
f(\lambda,0,0)=\sum_{\ell=0}^{r}\sum_{n=1}^{n_{\ell}}f_{\ell n}(\lambda)%
\mathsf{T}^{2\ell,n,0}
\end{equation}
with $f_{\ell n}(0)=0$ for $\ell>0$ and $1\leq n\leq n_{\ell}$. By %
\eqref{eq:11} we obtain
\begin{equation*}
f(\lambda,\varphi_{\mathbf{p}},\theta_{\mathbf{p}})=\sum_{\ell=0}^{r}
\sum_{n=1}^{n_{\ell}}f_{\ell n}(\lambda)\sum_{m=-2\ell}^{2\ell} U^{2\ell
g}_{m0}(\varphi_{\mathbf{p}},\theta_{\mathbf{p}})\mathsf{T}^{2\ell,n,m}.
\end{equation*}

Equation~\eqref{eq:12} takes the form
\begin{equation}  \label{eq:13}
\begin{aligned}
\langle\mathsf{T}(\mathbf{x}),\mathsf{T}(\mathbf{y})\rangle&=\frac{1}{2%
\sqrt{\pi}}
\sum_{\ell=0}^{r}\sum_{n=1}^{n_{\ell}}\sum_{m=-2\ell}^{2\ell}\int_{0}^{%
\infty}\int_{S^2} e^{\mathrm{i}(\mathbf{p},\mathbf{y}-\mathbf{x})}f_{\ell
n}(\lambda) \frac{1}{\sqrt{4\ell+1}}\\ &\quad\times
S^m_{2\ell}(\varphi_{\mathbf{p}},\theta_{\mathbf{p}})
\mathsf{T}^{2\ell,n,m}\sin\theta_{\mathbf{p}}\,\mathrm{d}\varphi_{%
\mathbf{p}}\, \mathrm{d}\theta_{\mathbf{p}}\,\mathrm{d}\nu(\lambda),
\end{aligned}
\end{equation}
where we used the relation
\begin{equation*}
U^{2\ell g}_{m0}(\varphi_{\mathbf{p}},\theta_{\mathbf{p}})=\sqrt{\frac{4\pi}{%
4\ell+1}} S^m_{2\ell}(\varphi_{\mathbf{p}},\theta_{\mathbf{p}}).
\end{equation*}

Substitute the \emph{Rayleigh expansion}
\begin{equation*}
\mathrm{e}^{\mathrm{i}(\mathbf{p},\mathbf{r})}=4\pi\sum^{\infty}_{\ell=0}
\sum^{\ell}_{m=-\ell}\mathrm{i}^{\ell}j_{\ell}(\|\mathbf{p}\|\cdot\|\mathbf{r%
}\|) S^m_{\ell}(\theta_{\mathbf{p}},\varphi_{\mathbf{p}}) S^m_{\ell}(\theta_{%
\mathbf{r}},\varphi_{\mathbf{r}})
\end{equation*}
into \eqref{eq:13}. We obtain
\begin{equation*}
\begin{aligned}
\langle\mathsf{T}(\mathbf{x}),\mathsf{T}(\mathbf{y})\rangle&=2\sqrt{\pi}
\sum_{\ell=0}^{r}\sum_{n=1}^{n_{\ell}}\sum_{m=-2\ell}^{2\ell}\int_{0}^{%
\infty} (-1)^{\ell}j_{2\ell}(\lambda\|\mathbf{r}\|)f_{\ell
n}(\lambda)\frac{1}{\sqrt{4\ell+1}}\\ &\quad\times
S^m_{\ell}(\varphi_{\mathbf{r}},\theta_{\mathbf{r}})
\mathsf{T}^{2\ell,n,m}\,\mathrm{d}\nu(\lambda), \end{aligned}
\end{equation*}
where $\mathbf{r}=\mathbf{y}-\mathbf{x}$. Returning back to the matrix
entries $U^{2\ell g}_{m0}(\varphi_{\mathbf{r}},\theta_{\mathbf{r}})$, we
have
\begin{equation}  \label{eq:14}
\langle\mathsf{T}(\mathbf{x}),\mathsf{T}(\mathbf{y})\rangle=
\int_{0}^{\infty}\sum_{\ell=0}^{r}(-1)^{\ell}j_{2\ell}(\lambda\|\mathbf{r}%
\|) \sum_{n=1}^{n_{\ell}}f_{\ell n}(\lambda) M^{2\ell,n}(\varphi_{\mathbf{r}%
},\theta_{\mathbf{r}})\,\mathrm{d}\nu(\lambda),
\end{equation}
where
\begin{equation*}
M^{2\ell,n}(\varphi_{\mathbf{r}},\theta_{\mathbf{r}})=\sum_{m=-2\ell}^{2%
\ell} U^{2\ell g}_{m0}(\varphi_{\mathbf{r}},\theta_{\mathbf{r}})\mathsf{T}%
^{2\ell,n,m}.
\end{equation*}

It is easy to check that the function $M^{2\ell,n}(\varphi_{\mathbf{r}%
},\theta_{\mathbf{r}})$ is a covariant of degree~$2\ell$ and of order~$2r$.
Therefore, the \emph{$M$-function} is a linear combination of basic
symmetric covariant tensors, or \emph{$L$-functions}:
\begin{equation*}
M^{2\ell,n}(\varphi_{\mathbf{r}},\theta_{\mathbf{r}})=\sum_{k=0}^{\ell}
\sum_{q=1}^{q_{kr}}c_{nkq}\frac{L^{2k,q}(\mathbf{y}-\mathbf{x})} {\|\mathbf{y%
}-\mathbf{x}\|^{2k}},
\end{equation*}
where $q_{kr}$ is the number of linearly independent symmetric covariant
tensors of degree~$2k$ and of order~$2r$. The right hand side is indeed a
polynomial in sines and cosines of the angles $\varphi_{\mathbf{r}}$ and $%
\theta_{\mathbf{r}}$. Equation~\eqref{eq:14} takes the form
\begin{equation*}
\begin{aligned}
\langle\mathsf{T}(\mathbf{x}),\mathsf{T}(\mathbf{y})\rangle&=
\int_{0}^{\infty}\sum_{\ell=0}^{r}(-1)^{\ell}j_{2\ell}(\lambda\|\mathbf{r}%
\|) \sum_{n=1}^{n_{\ell}}f_{\ell n}(\lambda)\\
&\quad\times \sum_{k=0}^{\ell}
\sum_{q=1}^{q_{kr}}c_{nkq}\frac{L^{2k,q}(\mathbf{y}-\mathbf{x})} {\|\mathbf{y%
}-\mathbf{x}\|^{2k}}\,\mathrm{d}\nu(\lambda).
\end{aligned}
\end{equation*}

Recall that $f_{\ell n}(\lambda)$ are measurable functions such that the
tensor \eqref{eq:15} lies in $\mathcal{C}_1$ for $\lambda=0$ and in $%
\mathcal{C}_0$ for $\lambda>0$. The final form of the two-point correlation
tensor of the random field $\mathsf{T}(\mathbf{x})$ is determined by
geometry of convex compacta $\mathcal{C}_0$ and $\mathcal{C}_1$. For
example, in the case of $r=1$ the set $\mathcal{C}_0$ is an interval (see %
\citet{MR3493458}), while $\mathcal{C}_1$ is a one-point set inside this
interval. The set $\mathcal{C}_0$ has two extreme points, and the
corresponding random field is a sum of two uncorrelated components given by
Equation~\eqref{eq:b1b2} below. The one-point set $\mathcal{C}_1$ lies in
the middle of the interval, the condition $\Phi_1(\{0\})=\Phi_2(\{0\})$
follows. In the case of $r=2$, the set of extreme points of the set $%
\mathcal{C}_0$ has three connected components: two one-point sets and an
ellipse, see \citet{MR3493458}, and the corresponding random field is a sum
of three uncorrelated components.

In general, the two-point correlation tensor of the field has the simplest
form when the set $\mathcal{C}_0$ is a simplex. We use this idea in Examples~%
\ref{ex:2components} and \ref{ex:5component} below.

\section{Examples of Mat\'{e}rn homogeneous and isotropic random fields}

\label{sec:examples}

\begin{example}
\label{ex:2components} Consider a centred homogeneous scalar isotropic
random field $T(\mathbf{x})$ on the space $\mathbb{R}^3$ with values in the
two-dimensional space $\mathbb{R}^2$. It is easy to see that both $\mathcal{C%
}_0$ and $\mathcal{C}_1$ are equal to the set of all symmetric
nonnegative-definite $2\times 2$ matrices with unit trace. Every such matrix
has the form
\begin{equation*}
\begin{pmatrix}
x & y \\
y & 1-x%
\end{pmatrix}%
\end{equation*}
with $x\in[0,1]$ and $y^2\leq x(1-x)$. Geometrically, $\mathcal{C}_0$ and $%
\mathcal{C}_1$ are the balls
\begin{equation*}
\left(x-\frac{1}{2}\right)^2+y^2=\frac{1}{4}.
\end{equation*}
Inscribe an equilateral triangle with vertices
\begin{equation*}
C^1=
\begin{pmatrix}
0 & 0 \\
0 & 1%
\end{pmatrix}
,\qquad C^{2,3}=\frac{1}{4}
\begin{pmatrix}
1 & \pm\sqrt{3} \\
\pm\sqrt{3} & 3%
\end{pmatrix}%
\end{equation*}
into the above ball. The function $f(\mathbf{p})$ takes the form
\begin{equation*}
f(\mathbf{p})=\sum_{m=1}^{3}a_m(\|\mathbf{p}\|)C^m,
\end{equation*}
where $a_m(\|\mathbf{p}\|)$ are the barycentric coordinates of the point $f(%
\mathbf{p})$ inside the triangle. The two-point correlation tensor of the
field takes the form
\begin{equation*}
\langle T(\mathbf{x}),T(\mathbf{y})\rangle=\sum_{m=1}^{3}\int_{0}^{\infty}
\frac{\sin(\lambda\|\mathbf{y}-\mathbf{x}\|)}{\lambda\|\mathbf{y}-\mathbf{x}%
\|} C^m\,\mathrm{d}\Phi_m(\lambda),
\end{equation*}
where $\mathrm{d}\Phi_m(\lambda)=a_m(\lambda)\mathrm{d}\nu(\lambda)$ are
three finite measures on $[0,\infty)$, and $\nu$ is the measure of Equation~%
\eqref{eq:12}. Define $\mathrm{d}\Phi_m(\lambda)$ as Mat\'{e}rn spectral
densities of Example~\ref{ex:2} (resp. dual Mat\'{e}rn spectral densities of
Example~\ref{ex:dual}). We obtain a scalar homogeneous and isotropic Mat\'{e}%
rn (resp. dual Mat\'{e}rn) random field.
\end{example}

\begin{example}
\label{ex:2component} Using \eqref{eq:Yaglom} and the well-known formulae
\begin{equation*}
j_0(t)=\frac{\sin t}{t},\qquad j_1(t)=\frac{\sin t}{t^2}-\frac{\cos t}{t},
\qquad j_2(t)=\left(\frac{3}{t^2}-1\right)\frac{\sin t}{t}-\frac{3\cos t}{t^2%
},
\end{equation*}
we write the two-point correlation tensor of rank~$1$ homogeneous and
isotropic random field in the form
\begin{equation*}
\langle\mathbf{v}(\mathbf{x}),\mathbf{v}(\mathbf{y})\rangle =B^{(1)}_{ij}(%
\mathbf{r})+B^{(2)}_{ij}(\mathbf{r}),
\end{equation*}
where $\mathbf{r}=\mathbf{y}-\mathbf{x}$, and
\begin{equation}  \label{eq:b1b2}
\begin{aligned}
B^{(1)}_{ij}(\mathbf{x},\mathbf{y})&=\int_{0}^{\infty}\left[\left(
-\frac{3\sin(\lambda\|\mathbf{r}\|)}{(\lambda\|\mathbf{r}\|)^3}
+\frac{\sin(\lambda\|\mathbf{r}\|)}{\lambda\|\mathbf{r}\|}
+\frac{3\cos(\lambda\|\mathbf{r}\|)}{(\lambda\|\mathbf{r}\|)^2}
\right)\frac{r_ir_j}{\|\mathbf{r}\|^2}\right.\\
&\quad+\left.\left(\frac{\sin(\lambda\|\mathbf{r}\|)}{(\lambda\|\mathbf{r}%
\|)^3}
-\frac{\cos(\lambda\|\mathbf{r}\|)}{(\lambda\|\mathbf{r}\|)^2}\right)%
\delta_{ij} \right]\,\mathrm{d}\Phi_1(\lambda),\\
B^{(2)}_{ij}(\mathbf{x},\mathbf{y})&=\int_{0}^{\infty}\left[\left(
\frac{3\sin(\lambda\|\mathbf{r}\|)}{(\lambda\|\mathbf{r}\|)^3}
-\frac{\sin(\lambda\|\mathbf{r}\|)}{\lambda\|\mathbf{r}\|}
-\frac{3\cos(\lambda\|\mathbf{r}\|)}{(\lambda\|\mathbf{r}\|)^2}
\right)\frac{r_ir_j}{\|\mathbf{r}\|^2}\right.\\
&\quad+\left.\left(\frac{\sin(\lambda\|\mathbf{r}\|)}{\lambda\|\mathbf{r}\|}
-\frac{\sin(\lambda\|\mathbf{r}\|)}{(\lambda\|\mathbf{r}\|)^3}
+\frac{\cos(\lambda\|\mathbf{r}\|)}{(\lambda\|\mathbf{r}\|)^2}\right)%
\delta_{ij} \right]\,\mathrm{d}\Phi_2(\lambda). \end{aligned}
\end{equation}

Now assume that the measures $\Phi_1$ and $\Phi_2$ are described by Mat\'{e}%
rn densities:
\begin{equation*}
\mathrm{d}\Phi_i(\lambda)=2\pi\lambda^2\frac{\sigma_i^{2}\Gamma \left( \nu_i
+\frac{3}{2}\right) a_i^{2\nu_i}}{2\pi ^{3/2}\left(
a_i^{2}+\lambda^{2}\right) ^{\nu_i +\frac{3}{2}}},\qquad i=1,2.
\end{equation*}
It is possible to substitute these densities to \eqref{eq:b1b2} and
calculate the integrals using \citet[Equation~2.5.9.1]{MR874986}. We obtain
rather long expressions that include the generalised hypergeometric function
${}_1F_2$.

The situation is different for the dual model:
\begin{equation*}
\mathrm{d}\Phi_i(\lambda)=\frac{1}{(2\pi)^22^{\nu_i-1}\Gamma(\nu_i+3/2)}
\lambda^{\nu+2}K_{\nu}(\lambda).
\end{equation*}
Using \citet[Equations 2.16.14.3, 2.16.14.4]{MR950173}, we obtain
\begin{equation*}
\begin{aligned}
B^{(1)}_{ij}(\mathbf{x},\mathbf{y})&=C_1\left(-\frac{3\pi\Gamma(2\nu_1)}{4\|%
\mathbf{r}\|^3 (1+\|\mathbf{r}\|^2)^{\nu_1/2}}
\left[P^{-\nu_1}_{\nu_1-1}\left(
\frac{\|\mathbf{r}\|}{\sqrt{1+\|\mathbf{r}\|^2}}\right)\right.\right.\\
&\quad\left.\left.-P^{-\nu_1}_{\nu_1-1}\left(-\frac{\|\mathbf{r}\|}
{\sqrt{1+\|\mathbf{r}\|^2}}\right)\right]+\frac{2^{\nu_1}\sqrt{\pi}\Gamma(%
\nu_1+3/2)} {(1+\|\mathbf{r}\|^2)^{\nu_1+3/2}}\right.\\
&\quad+\left.\frac{3\cdot 2^{\nu_1-1}\sqrt{\pi}\Gamma(\nu_1+1/2)}
{(1+\|\mathbf{r}\|^2)^{\nu_1+1/2}}\right)\frac{r_ir_j}{\|\mathbf{r}\|^2}\\
&\quad+C_1\left(\frac{\pi\Gamma(2\nu_1)}{4\|\mathbf{r}\|^3
(1+\|\mathbf{r}\|^2)^{\nu_1/2}}\left[P^{-\nu_1}_{\nu_1-1}\left(
\frac{\|\mathbf{r}\|}{\sqrt{1+\|\mathbf{r}\|^2}}\right)\right.\right.\\
&\quad\left.\left.-P^{-\nu_1}_{\nu_1-1}\left(-\frac{\|\mathbf{r}\|}
{\sqrt{1+\|\mathbf{r}\|^2}}\right)\right]-\frac{2^{\nu_1}\sqrt{\pi}\Gamma(%
\nu_1+1/2)} {(1+\|\mathbf{r}\|^2)^{\nu_1+3/2}}\right)\delta_{ij},\\
\end{aligned}
\end{equation*}
and
\begin{equation*}
\begin{aligned}
B^{(2)}_{ij}(\mathbf{x},\mathbf{y})&=C_2\left(\frac{3\pi\Gamma(2\nu_2)}{4\|%
\mathbf{r}\|^3 (1+\|\mathbf{r}\|^2)^{\nu_2/2}}
\left[P^{-\nu_2}_{\nu_2-1}\left(
\frac{\|\mathbf{r}\|}{\sqrt{1+\|\mathbf{r}\|^2}}\right)\right.\right.\\
&\quad\left.\left.-P^{-\nu_2}_{\nu_2-1}\left(-\frac{\|\mathbf{r}\|}
{\sqrt{1+\|\mathbf{r}\|^2}}\right)\right]-\frac{2^{\nu_2}\sqrt{\pi}\Gamma(%
\nu_2+3/2)} {(1+\|\mathbf{r}\|^2)^{\nu_2+3/2}}\right.\\
&\quad-\left.\frac{3\cdot 2^{\nu_2-1}\sqrt{\pi}\Gamma(\nu_2+1/2)}
{(1+\|\mathbf{r}\|^2)^{\nu_2+1/2}}\right)\frac{r_ir_j}{\|\mathbf{r}\|^2}\\
&\quad+C_2\left(\frac{\sqrt{\pi}\Gamma(\nu_2+3/2)2^{\nu_2}}{\|\mathbf{r}\|
(1+\|\mathbf{r}\|^2)^{\nu_2+3/2}}-\frac{\pi\Gamma(2\nu_2)}{4\|\mathbf{r}\|^3
(1+\|\mathbf{r}\|^2)^{\nu_2/2}}\right.\\
&\quad\times\left[P^{-\nu_2}_{\nu_2-1}\left(
\frac{\|\mathbf{r}\|}{\sqrt{1+\|\mathbf{r}\|^2}}\right)-P^{-\nu_2}_{\nu_2-1}
\left(-\frac{\|\mathbf{r}\|}{\sqrt{1+\|\mathbf{r}\|^2}}\right)\right]\\
&\quad-\left.\frac{2^{\nu_2-1}\sqrt{\pi}\Gamma(\nu_2+1/2)}
{\|\mathbf{r}\|^2(1+\|\mathbf{r}\|^2)^{\nu_2+1/2}}\right)\delta_{ij},\\
\end{aligned}
\end{equation*}
where
\begin{equation*}
C_i=\frac{1}{(2\pi)^22^{\nu_i-1}\Gamma(\nu_i+3/2)}.
\end{equation*}
\end{example}

\begin{example}
\label{ex:5component} Consider the case when $r=2$ and $U(g)=\mathsf{S}^2(g)$%
. In order to write down symmetric rank~$4$ tensors in a compressed matrix
form, consider an orthogonal operator $\tau$ acting from $\mathsf{S}^2(%
\mathsf{S}^2(\mathbb{R}^3))$ to $\mathsf{S}^2(\mathbb{R}^6)$ as follows:
\begin{equation*}
\tau f_{ijkl}=\left(
\begin{smallmatrix}
f_{-1-1-1-1} & f_{-1-100} & f_{-1-111} & \sqrt{2}f_{-1-1-10} & \sqrt{2}%
f_{-1-101} & \sqrt{2}f_{-1-11-1} \\
f_{00-1-1} & f_{0000} & f_{0011} & \sqrt{2}f_{00-10} & \sqrt{2}f_{0001} &
\sqrt{2}f_{001-1} \\
f_{11-1-1} & f_{1100} & f_{1111} & \sqrt{2}f_{11-10} & \sqrt{2}f_{1101} &
\sqrt{2}f_{111-1} \\
\sqrt{2}f_{-10-1-1} & \sqrt{2}f_{-1000} & \sqrt{2}f_{-1011} & 2f_{-10-10} &
2f_{-1001} & 2f_{-101-1} \\
\sqrt{2}f_{01-1-1} & \sqrt{2}f_{0100} & \sqrt{2}f_{0111} & 2f_{01-10} &
2f_{0101} & 2f_{011-1} \\
\sqrt{2}f_{1-1-1-1} & \sqrt{2}f_{1-100} & \sqrt{2}f_{1-111} & 2f_{1-1-10} &
2f_{1-101} & 2f_{1-11-1}%
\end{smallmatrix}
\right),
\end{equation*}
see \cite[Equation~(44)]{MR1816224}. It is possible to prove the following.
The matrix $\tau f_{ijkl}(\mathbf{0})$ lies in the interval $\mathcal{C}_1$
with extreme points $C^1$ and $C^2$, where the nonzero elements of the
symmetric matrix $C^1$ lying on and over the main diagonal are as follows:
\begin{equation*}
C^1_{11}=C^1_{12}=C^1_{13}=C^1_{22}=C^1_{23}=C^1_{33}=\frac{1}{3},
\end{equation*}
while those of the matrix $C^2$ are
\begin{equation*}
\begin{aligned}
C^2_{11}&=C^2_{22}=C^2_{33}=\frac{2}{15},\qquad C^2_{44}=C^2_{55}=C^2_{66}=%
\frac{1}{5},\\
C^2_{12}&=C^2_{13}=C^2_{23}=-\frac{1}{15}.
\end{aligned}
\end{equation*}
The matrix $\tau f_{ijkl}(\lambda,0,0)$ with $\lambda>0$ lies in the convex
compact set $\mathcal{C}_0$. The set of extreme points of $\mathcal{C}_0$
contains three connected components. The first component is the one-point
set $\{D^1\}$ with
\begin{equation*}
D^1_{44}=D^1_{66}=\frac{1}{2}.
\end{equation*}
The second component is the one-point set $\{D^2\}$ with
\begin{equation*}
D^2_{11}=D^2_{33}=\frac{1}{4},\qquad D^2_{55}=\frac{1}{2},\qquad D^2_{13}=-%
\frac{1}{4}.
\end{equation*}
The third component is the ellipse $\{\,D^{\theta}\colon
0\leq\theta<2\pi\,\} $ with
\begin{equation*}
\begin{aligned}
D^{\theta}_{11}&=D^{\theta}_{33}=D^{\theta}_{13}=\frac{1}{2}%
\sin^2(\theta/2),\qquad D^{\theta}_{22}=\cos^2(\theta/2),\\
D^{\theta}_{12}&=D^{\theta}_{23}=\frac{1}{2\sqrt{2}}\sin(\theta).
\end{aligned}
\end{equation*}

Choose three points $D^3$, $D^4$, $D^5$ lying on the above ellipse. If we
allow the matrix $\tau f_{ijkl}(\lambda,0,0)$ with $\lambda>0$ to take
values in the simplex with vertices $D^i$, $1\leq i\leq 5$, then the
two-point correlation tensor of the random field $\varepsilon(\mathbf{x})$
is the sum of five integrals. The more the four-dimensional Lebesgue measure
of the simplex in comparison with that of $\mathcal{C}_0$, the wider class
of random fields is described.

Note that the simplex should contain the set $\mathcal{C}_1$. The matrix $%
C^1 $ lies on the ellipse and corresponds to the value of $\theta=2\arcsin(%
\sqrt{2/3})$. It follows that one of the above points, say $D^3$, must be
equal to $C^1$. If we choose $D^4$ to correspond to the value of $%
\theta=2(\pi-\arcsin(\sqrt{2/3}))$, that is,
\begin{equation*}
D^4_{11}=D^4_{33}=D^4_{13}=\frac{1}{6},\qquad D^4_{22}=\frac{2}{3},\qquad
D^4_{12}=D^4_{23}=-\frac{1}{3},
\end{equation*}
then
\begin{equation*}
C^2=\frac{2}{5}(D^1+D^2)+\frac{1}{5}D^4,
\end{equation*}
and $C^2$ lies in the simplex. Finally, choose $D_5$ to correspond to the
value of $\theta=\pi$, that is
\begin{equation*}
D^5_{11}=D^5_{33}=D^5_{13}=\frac{1}{2}.
\end{equation*}
The constructed simplex is not the one with maximal possible Lebesgue
measure, but the coefficients in formulas are simple.

\begin{theorem}
Let $\varepsilon(\mathbf{x})$ be a random field that describes the stress
tensor of a deformable body. The following conditions are equivalent.

\begin{enumerate}
\item The matrix $\tau f_{ijkl}(\lambda,0,0)$ with $\lambda>0$ takes values
in the simplex described above.

\item The correlation tensor of the field has the spectral expansion
\begin{equation*}
\langle\varepsilon(\mathbf{x}),\varepsilon(\mathbf{y})\rangle=
\sum^5_{n=1}\int^{\infty}_0\sum^5_{q=1} \tilde{N}_{nq}(\lambda,\|\mathbf{r}%
\|)L^q_{ijkl}(\mathbf{r})\,\mathrm{d}\Phi_n(\lambda),
\end{equation*}
where the non-zero functions $\tilde{N}_{nq}(\lambda,r)$ are given in Table~%
\emph{\ref{tab:3}}, and where $\Phi_n(\lambda)$ are five finite measures on $%
[0,\infty)$ with
\begin{equation*}
\Phi_1(\{0\})=\Phi_2(\{0\})=2\Phi_4(\{0\}),\qquad \Phi_5(\{0\})=0.
\end{equation*}
\end{enumerate}
\end{theorem}

Assume that all measures $\Phi_n$ are absolutely continuous and their
densities are either the Mat\'{e}rn or the dual Mat\'{e}rn densities. The
two-point correlation tensors of the corresponding random fields can be
calculated in exactly the same way as in Example~\ref{ex:2component}.

\begin{table}[tbp]
\caption{The functions $\tilde{N}_{nq}(\protect\lambda,r)$}
\label{tab:3}%
\begin{tabular}{|l|l|l|}
\hline
$n$ & $q$ & $N_{nq}(\lambda,r)$ \\ \hline
1 & 1 & $-\frac{1}{15}j_0(\lambda r)-\frac{2}{21}j_2(\lambda r)-\frac{1}{35}%
j_4(\lambda r)$ \\
1 & 2 & $\frac{1}{10}j_0(\lambda r)+\frac{1}{14}j_2(\lambda r)-\frac{1}{35}%
j_4(\lambda r)$ \\
1 & 3 & $-\frac{3}{28}j_2(\lambda r)+\frac{1}{7}j_4(\lambda r)$ \\
1 & 4 & $\frac{1}{7}j_2(\lambda r)+\frac{1}{7}j_4(\lambda r)$ \\
1 & 5 & $-j_4(\lambda r)$ \\
2 & 1 & $-\frac{1}{15}j_0(\lambda r)+\frac{4}{21}j_2(\lambda r)+\frac{1}{140}%
j_4(\lambda r)$ \\
2 & 2 & $\frac{1}{10}j_0(\lambda r)-\frac{1}{7}j_2(\lambda r)+\frac{1}{140}%
j_4(\lambda r)$ \\
2 & 3 & $\frac{3}{14}j_2(\lambda r)-\frac{1}{28}j_4(\lambda r)$ \\
2 & 4 & $-\frac{2}{7}j_2(\lambda r)-\frac{1}{28}j_4(\lambda r)$ \\
2 & 5 & $\frac{1}{4}j_4(\lambda r)$ \\
3 & 1 & $\frac{1}{3}j_0(\lambda r)$ \\
4 & 1 & $-\frac{1}{135}j_0(\lambda r)-\frac{4}{21}j_2(\lambda r)+\frac{3}{70}%
j_4(\lambda r)$ \\
4 & 2 & $\frac{1}{90}j_0(\lambda r)+\frac{1}{7}j_2(\lambda r)+\frac{3}{70}%
j_4(\lambda r)$ \\
4 & 3 & $-\frac{3}{14}j_2(\lambda r)-\frac{3}{14}j_4(\lambda r)$ \\
4 & 4 & $\frac{2}{7}j_2(\lambda r)-\frac{3}{14}j_4(\lambda r))$ \\
4 & 5 & $\frac{3}{2}j_4(\lambda r)$ \\
5 & 1 & $\frac{1}{5}j_0(\lambda r)-\frac{2}{7}j_2(\lambda r)+\frac{1}{70}%
j_4(\lambda r)$ \\
5 & 2 & $\frac{1}{30}j_0(\lambda r)+\frac{2}{21}j_2(\lambda r)+\frac{1}{70}%
j_4(\lambda r)$ \\
5 & 3 & $\frac{1}{14}j_2(\lambda r)-\frac{1}{14}j_4(\lambda r)$ \\
5 & 4 & $\frac{5}{21}j_2(\lambda r)-\frac{1}{14}j_4(\lambda r)$ \\
5 & 5 & $\frac{1}{2}j_4(\lambda r)$ \\ \hline
\end{tabular}%
\end{table}

Introduce the following notation:
\begin{equation*}
\begin{aligned}
\mathsf{T}^{0,1}_{ijkl}&=\frac{1}{3}\delta_{ij}\delta_{kl},\\
\mathsf{T}^{0,2}_{ijkl}&=\frac{1}{\sqrt{5}}%
\sum_{n=-2}^{2}g^{n[i,j]}_{2[1,1]} g^{n[k,l]}_{2[1,1]},\\
\mathsf{T}^{2,1,m}_{ijkl}&=\frac{1}{\sqrt{6}}(\delta_{ij}g^{m[k,l]}_{2[1,1]}
+\delta_{kl}g^{m[i,j]}_{2[1,1]}),\qquad -2\leq m\leq 2,\\
\mathsf{T}^{2,2,m}_{ijkl}&=%
\sum_{n,q=-2}^{2}g^{m[n,q]}_{2[2,2]}g^{n[i,j]}_{2[1,1]}
g^{q[k,l]}_{2[1,1]},\qquad -2\leq m\leq 2,\\
\mathsf{T}^{4,1,m}_{ijkl}&=%
\sum_{n,q=-4}^{4}g^{m[n,q]}_{4[2,2]}g^{n[i,j]}_{2[1,1]}
g^{q[k,l]}_{2[1,1]},\qquad -4\leq m\leq 4. \end{aligned}
\end{equation*}
Consider the five nonnegative-definite matrices $A^n$, $1\leq n\leq 5$, with
the following matrix entries:
\begin{equation*}
\begin{aligned}
a^{\ell''m''kl,1}_{\ell'm'ij}&=\sqrt{(2\ell'+1)(2\ell''+1)}\left(\frac{1}{%
\sqrt{5}}
\mathsf{T}^{0,2}_{ijkl}g^{0[m',m'']}_{0[\ell',\ell'']}g^{0[0,0]}_{0[\ell',%
\ell'']}\right.\\
&\quad-\left.\frac{1}{5\sqrt{14}}\sum_{m=-2}^{2}\mathsf{T}^{2,2,m}_{ijkl}
g^{m[m',m'']}_{2[\ell',\ell'']}g^{0[0,0]}_{2[\ell',\ell'']}
-\frac{2\sqrt{2}}{9\sqrt{35}}\sum_{m=-4}^{4}\mathsf{T}^{4,1,m}_{ijkl}
g^{m[m',m'']}_{4[\ell',\ell'']}g^{0[0,0]}_{4[\ell',\ell'']}\right),\\
a^{\ell''m''kl,2}_{\ell'm'ij}&=\sqrt{(2\ell'+1)(2\ell''+1)}\left(\frac{1}{%
\sqrt{5}}
\mathsf{T}^{0,2}_{ijkl}g^{0[m',m'']}_{0[\ell',\ell'']}g^{0[0,0]}_{0[\ell',%
\ell'']}\right.\\
&\quad+\left.\frac{\sqrt{2}}{5\sqrt{7}}\sum_{m=-2}^{2}%
\mathsf{T}^{2,2,m}_{ijkl}
g^{m[m',m'']}_{2[\ell',\ell'']}g^{0[0,0]}_{2[\ell',\ell'']}
+\frac{1}{9\sqrt{70}}\sum_{m=-4}^{4}\mathsf{T}^{4,1,m}_{ijkl}
g^{m[m',m'']}_{4[\ell',\ell'']}g^{0[0,0]}_{4[\ell',\ell'']}\right),\\
a^{\ell''m''kl,3}_{\ell'm'ij}&=\sqrt{(2\ell'+1)(2\ell''+1)}
\mathsf{T}^{0,1}_{ijkl}g^{0[m',m'']}_{0[\ell',\ell'']}g^{0[0,0]}_{0[\ell',%
\ell'']},\\
a^{\ell''m''kl,4}_{\ell'm'ij}&=\sqrt{(2\ell'+1)(2\ell''+1)}\left(\frac{1}{9%
\sqrt{5}}
\mathsf{T}^{0,2}_{ijkl}g^{0[m',m'']}_{0[\ell',\ell'']}g^{0[0,0]}_{0[\ell',%
\ell'']}\right.\\
&\quad-\left.\frac{\sqrt{2}}{5\sqrt{7}}\sum_{m=-2}^{2}%
\mathsf{T}^{2,2,m}_{ijkl}
g^{m[m',m'']}_{2[\ell',\ell'']}g^{0[0,0]}_{2[\ell',\ell'']}
+\frac{\sqrt{2}}{3\sqrt{35}}\sum_{m=-4}^{4}\mathsf{T}^{4,1,m}_{ijkl}
g^{m[m',m'']}_{4[\ell',\ell'']}g^{0[0,0]}_{4[\ell',\ell'']}\right),\\
a^{\ell''m''kl,5}_{\ell'm'ij}&=\sqrt{(2\ell'+1)(2\ell''+1)}\left(\left(%
\frac{2}{3} \mathsf{T}^{0,1}_{ijkl}+\frac{1}{3\sqrt{5}}
\mathsf{T}^{0,2}_{ijkl}\right)g^{0[m',m'']}_{0[\ell',\ell'']}g^{0[0,0]}_{0[%
\ell',\ell'']}\right.\\
&\quad+\left(\frac{2}{9}\sum_{m=-2}^{2}\mathsf{T}^{2,1,m}_{ijkl}
-\frac{\sqrt{2}}{9\sqrt{7}}\sum_{m=-2}^{2}\mathsf{T}^{2,2,m}_{ijkl}
\right)g^{m[m',m'']}_{2[\ell',\ell'']}g^{0[0,0]}_{2[\ell',\ell'']}\\
&\quad+\left.\frac{\sqrt{2}}{9\sqrt{35}}\sum_{m=-4}^{4}%
\mathsf{T}^{4,1,m}_{ijkl}
g^{m[m',m'']}_{4[\ell',\ell'']}g^{0[0,0]}_{4[\ell',\ell'']}\right),
\end{aligned}
\end{equation*}
and let $L^n$ be infinite lower triangular matrices from Cholesky
factorisation of the matrices $A^n$.

\begin{theorem}
The following conditions are equivalent.

\begin{enumerate}
\item The matrix $\tau f_{ij\ell m}(\lambda,0,0)$ with $\lambda>0$ takes
values in the simplex described above.

\item The field $\varepsilon(\mathbf{x})$ has the form
\begin{equation*}
\varepsilon_{ij}(\rho,\theta,\varphi)=C\delta_{ij}+2\sqrt{\pi}%
\sum_{n=1}^{5}\sum_{\ell=0}^{\infty}
\sum_{m=-\ell}^{\ell}\int_{0}^{\infty}j_{\ell}(\lambda\rho)\,\mathrm{d}%
Z^{n^{\prime}}_{\ell mij}(\lambda)S^m_{\ell}(\theta,\varphi),
\end{equation*}
where
\begin{equation*}
Z^{n^{\prime}}_{\ell
mij}(A)=\sum_{(\ell^{\prime},m^{\prime},k,l)\leq(\ell,m,i,j)}Z^n_{\ell^{%
\prime}m^{\prime}kl}(A),
\end{equation*}
and where $Z^n_{\ell^{\prime}m^{\prime}kl}$ is the sequence of uncorrelated
scattered random measures on $[0,\infty)$ with control measures $\Phi_n$.
\end{enumerate}
\end{theorem}

The idea of proof is as follows. Write down the Rayleigh expansion for $%
\mathrm{e}^{\mathrm{i}(\mathbf{p},\mathbf{x})}$ and for $\mathrm{e}^{-%
\mathrm{i}(\mathbf{p},\mathbf{y})}$ separately,substitute both expansions
into \eqref{eq:13} and use the following result, known as the \emph{Gaunt
integral}:
\begin{equation*}
\begin{aligned}
\int_{S^2}S^{m_1}_{\ell_1}(\theta,\varphi)S^{m_2}_{\ell_2}(\theta,\varphi)
S^{m_3}_{\ell_3}(\theta,\varphi)\sin\theta\,\mathrm{d}\varphi\, \mathrm{d}%
\theta&=\sqrt{\frac{(2\ell_1+1)(2\ell_2+1)}{4\pi(2\ell_3+1)}}\\
&\quad\times
g^{m_3[m_1,m_2]}_{\ell_3[\ell_1,\ell_2]}g^{0[0,0]}_{\ell_3[\ell_1,\ell_2]}.
\end{aligned}
\end{equation*}
This theorem can be proved exactly in the same way, as its complex
counterpart, see, for example, \citet{MR2840154}. Then apply Karhunen's
theorem, see \citet{MR0023013}.
\end{example}

\bigskip

\textbf{Acknowledgements}. Nikolai N. Leonenko was supported in part by
projects MTM2012-32674 (co-funded by European Regional Development Funds),
and MTM2015--71839--P, MINECO, Spain. This research was also supported under
Australian Research Council's Discovery Projects funding scheme (project
number DP160101366), and under Cardiff Incoming Visiting Fellowship Scheme
and International Collaboration Seedcorn Fund.

Anatoliy Malyarenko is grateful to Professor Martin Ostoja-Starzewski for useful him to probabilistic models of continuum physics and fruitful discussions.

\bibliography{\jobname}

\appendix

\section{Tensors}

\label{ap:tensors}

There are several equivalent definitions of tensors. Surprisingly, the most
abstract of them is useful in the theory of random fields.

Let $r$ be a nonnegative integer, and let $V_1$, \dots, $V_r$ be linear
spaces over the same field $\mathbb{K}$. When $r=0$, define the tensor
product of the empty family of spaces as $\mathbb{K}^1$, the one-dimensional
linear space over $\mathbb{K}$.

\begin{theorem}[The universal mapping property]
There exist a unique linear space $V_1\otimes\cdots\otimes V_r$ and a unique
linear operator $\tau\colon V_1\times V_2\times\cdots\times V_r\to
V_1\otimes\cdots\otimes V_r$ that satisfy the \emph{universal mapping
property}: for any linear space $W$ and for any multilinear map $\beta\colon
V_1\times V_2\times\cdots\times V_r\to W$, there exists a unique \emph{linear%
} operator $B\colon V_1\otimes\cdots\otimes V_r\to X$ such that $%
\beta=B\circ\tau$:
\begin{equation*}
\xymatrix{V_1\times V_2\times\cdots\times V_r\ar[dr]_{\beta}\ar[r]^{\tau} &
V_1\otimes\cdots\otimes V_r\ar[d]^B\\ & W}
\end{equation*}
\end{theorem}

In other words: \emph{the construction of the tensor product of linear
spaces reduces the study of multilinear mappings to the study of linear ones}%
.

The tensor product $\mathbf{v}_1\otimes\cdots\otimes\mathbf{v}_r$ of the
vectors $\mathbf{v}_i\in V_i$, $1\leq i\leq r$, is defined by
\begin{equation*}
\mathbf{v}_1\otimes\cdots\otimes\mathbf{v}_r=\tau(\mathbf{v}_1,\dots,\mathbf{%
v}_r).
\end{equation*}

Let $V_1$, \dots, $V_r$, $W_1$, \dots, $W_r$ be finite-dimensional linear
spaces, and let $A_i\in L(V_i,W_i)$ for $1\leq i\leq r$. The tensor product
of linear operators, $A_1\otimes\cdots\otimes A_r$, is a unique element of
the space $L(V_1\otimes\cdots\otimes V_r,W_1\otimes\cdots\otimes W_r)$ such
that
\begin{equation*}
(A_1\otimes\cdots\otimes A_r)(\mathbf{v}_1\otimes\cdots\otimes\mathbf{v}%
_r):=A_1(\mathbf{v}_1)\otimes \cdots\otimes A_r(\mathbf{v}_r),\qquad\mathbf{v%
}_i\in V_i.
\end{equation*}

If all the spaces $V_i$, $1\leq i\leq r$, are copies of the same space $V$,
then we write $V^{\otimes r}$ for the $r$-fold tensor product of $V$ with
itself, and $\mathbf{v}^{\otimes r}$ for the tensor product of $r$ copies of
a vector $\mathbf{v}\in V$. Similarly, for $A\in L(V,V)$ we write $%
A^{\otimes r}$ for the $r$-fold tensor product of $A$ with itself. Note that
$A^{\otimes 0}$ is the identity operator in the space $\mathbb{K}^1$.

\section{Group representations}

Let $G$ be a topological group. A \emph{finite-dimensional representation}
of $G$ is a pair $(\rho,V)$, where $V$ is a finite-dimensional linear space,
and $\rho\colon G\to\mathrm{GL}(V)$ is a continuous group homomorphism. Here
$\mathrm{GL}(V)$ is the \emph{general linear group of order} $n$, or the
group of all invertible $n\times n$ matrices. In what follows, we omit the
word ``finite-dimensional'' unless infinite-dimensional representations are
under consideration.

In a coordinate form, a representation of $G$ is a continuous group
homomorphism $\rho\colon G\to\mathrm{GL}(n,\mathbb{K})$ and the space $%
\mathbb{K}^n$.

Let $W\subseteq V$ be a linear subspace of the space~$V$. $W$ is called an
\emph{invariant subspace} of the representation $(\rho,V)$ if $\rho(g)%
\mathbf{w}\in W$ for all $g\in G$ and $\mathbf{w}\in W$. The restriction of $%
\rho$ to $W$ is then a representation $(\sigma,W)$ of $G$. Formula
\begin{equation*}
\tau(g)(\mathbf{v}+W):=\rho(g)\mathbf{v}+W
\end{equation*}
defines a representation $(\tau,V/W)$ of $G$ in the quotient space $V/W$.

In a coordinate form, take a basis for $W$ and complete it to a basis for $V$%
. The matrix of $\rho(g)$ relative to the above basis is
\begin{equation}  \label{eq:matrixform}
\rho(g)=
\begin{pmatrix}
\sigma(g) & * \\
0 & \tau(g)%
\end{pmatrix}
.
\end{equation}

Let $(\rho,V)$ and $(\tau,W)$ be representations of $G$. An operator $A\in
L(V,W)$ is called an \emph{intertwining operator} if
\begin{equation}  \label{eq:intertwining}
\tau(g)A=A\rho(g),\qquad g\in G.
\end{equation}
The intertwining operators form a linear space $L_G(V,W)$ over $\mathbb{F}$.

The representations $(\rho,V)$ and $(\tau,W)$ are called \emph{equivalent}
if the space $L_G(V,W)$ contains an invertible operator. Let $A$ be such an
operator. Multiply \eqref{eq:intertwining} by $A^{-1}$ from the right. We
obtain
\begin{equation*}
\tau(g)=A\rho(g)A^{-1},\qquad g\in G.
\end{equation*}
In a coordinate form, $\tau(g)$ and $\rho(g)$ are matrices of the same
presentation, written in two different bases, and $A$ is the transition
matrix between the bases.

A representation $(\rho,V)$ with $V\neq\{\mathbf{0}\}$ is called \emph{%
reducible} if if there exists an invariant subspace $W\notin\{\{\mathbf{0}%
\},V\}$. In a coordinate form, all blocks of the matrix \eqref{eq:matrixform}
are nonempty. Otherwise, the representation is called \emph{irreducible}.

\begin{example}
Let $G=\mathrm{O}(3)$. The mapping $g\mapsto g^{\otimes r}$ is a
representation of the group~$G$ in the space $(\mathbb{R}^3)^{\otimes r}$.
When $r=0$, this representation is called \emph{trivial}, when $r=1$, it is
called \emph{defining}. When $r\geq 2$, this representation is reducible.
\end{example}

From now on we suppose that the topological group $G$ is compact. There
exists an inner product $(\boldsymbol{\cdot},\boldsymbol{\cdot})$ on $V$
such that
\begin{equation*}
(\rho(g)\mathbf{v},\rho(g)\mathbf{w})=(\mathbf{v},\mathbf{w}),\qquad\mathbf{v%
}, \mathbf{w}\in V.
\end{equation*}
In a coordinate form, we can choose an orthonormal basis in $V$. If $V$ is a
complex linear space, then the representation $(\rho,V)$ takes values in $%
\mathrm{U}(n)$, the group of $n\times n$ unitary matrices, and we speak of a
\emph{unitary representation} If $V$ is a real linear space, then the
representation $(\rho,V)$ takes values in $\mathrm{O}(n)$, and we speak of
an \emph{orthogonal representation}.

Let $(\pi,V)$ and $(\rho,W)$ be representations of $G$. The \emph{direct sum
of representations} is the representation $(\pi\oplus\rho,V\oplus W)$ acting
by
\begin{equation*}
(\pi\oplus\rho)(g)(\mathbf{v}\oplus\mathbf{w}):=\pi(g)\mathbf{v}%
\oplus\rho(g) \mathbf{w},\qquad g\in G,\quad\mathbf{v}\in V,\quad\mathbf{w}%
\in W.
\end{equation*}
In a coordinate form, we have
\begin{equation}  \label{eq:reducible}
\pi\oplus\rho(g)=
\begin{pmatrix}
\pi(g) & 0 \\
0 & \rho(g)%
\end{pmatrix}
.
\end{equation}

Consider the action $\pi\otimes\rho$ of the group $G$ on the set of tensor
products $\mathbf{v}\otimes\mathbf{w}$ defined by
\begin{equation*}
(\pi\otimes\rho)(g)(\mathbf{v}\otimes\mathbf{w}):=\pi(g)\mathbf{v}%
\otimes\rho(g) \mathbf{w},\qquad g\in G,\quad\mathbf{v}\in V,\quad\mathbf{w}%
\in W.
\end{equation*}
This action may be extended by linearity to the \emph{tensor product of
representations} $(\pi\otimes\rho,V\otimes W)$. In a coordinate form, $%
(\pi\otimes\rho)(g)$ is a rank $4$ tensor with components
\begin{equation*}
\mathsf{T}_{ijkl}(g)=\pi_{ij}(g)\rho_{kl}(g),\qquad g\in G.
\end{equation*}

A representation $(\sigma,V)$ of a group $G$ is called \emph{completely
reducible} if for every invariant subspace $W\subset V$ there exists an
invariant subspace $U\subset V$ such that $V=W\oplus U$. In a coordinate
form, any basis $\{\mathbf{w}_1,\dots,\mathbf{w}_p\}$ for $W$ can be
completed to a basis $\{\mathbf{w}_1,\dots,\mathbf{w}_p,\mathbf{u}_1,\dots,%
\mathbf{u}_q\}$ for $V$ such that the span of the vectors $\mathbf{u}_1$%
,\dots,$\mathbf{u}_q$ is invariant. The matrix $\sigma(g)$ in the above
basis has the form \eqref{eq:reducible}. Any representation of a compact
group is completely reducible.

Let $(\rho,V)$ be an irreducible representation of a group $G$. Denote by $%
[\rho]$ the equivalence class of all representations of $G$ equivalent to $%
(\rho,V)$ and by $\hat{G}$ the set of all equivalence classes of irreducible
representations of $G$. For any finite-dimensional representation $%
(\sigma,V) $ of $G$, there exists finitely many equivalence classes $[\rho_1]
$, \dots, $[\rho_k]\in\hat{G}$ and uniquely determined positive integers $m_1
$, \dots, $m_k$ such that $(\sigma,V)$ is equivalent to the direct sum of $%
m_1$ copies of the representation $(\rho_1,V_1)$, \dots, $m_k$ copies of the
representation $(\rho_k,V_k)$. The direct sum $m_iV_i$ of $m_i$ copies of
the linear space $V_i$ is called the \emph{isotypic subspace} of the space $%
V $ that corresponds to the representation $(\rho_i,V_i)$. The numbers $m_i$
are called the \emph{multiplicities} of the irreducible representation $%
(\rho_i,V_i)$ in $(\sigma,V)$. The decompositions $V=\sum m_iV_i$ and $%
\sigma=\sum m_i\rho_i$ are called the \emph{isotypic decompositions}.

Assume that a compact group $G$ is \emph{easy reducible}. This means that
for any three irreducible representation $(\rho,V)$, $(\sigma,W)$, and $%
(\tau,U)$ of $G$ the multiplicity $m_{\tau}$ of $\tau$ in $\rho\otimes\sigma$
is equal to either $0$ or $1$. For example, the group $\mathrm{O}(3)$ is
easy reducible. Assume $m_{\tau}=1$. Let $\{\,\mathbf{e}^{\rho}_i\colon
1\leq i\leq\dim\rho\,\}$ be an orthonormal basis in $V$, and similarly for $%
\sigma$ and $\tau$. There are two natural bases in the space $V\otimes W$.
The \emph{coupled basis} is
\begin{equation*}
\{\,\mathbf{e}^{\rho}_i\otimes\mathbf{e}^{\sigma}_j\colon 1\leq
i\leq\dim\rho,1\leq j\leq\dim\sigma\,\}.
\end{equation*}
The \emph{uncoupled basis} is
\begin{equation*}
\{\,\mathbf{e}^{\tau}_k\colon m_{\tau}=1,1\leq k\leq\dim\tau\,\}.
\end{equation*}

In a coordinate form, the elements of the space $V\otimes W$ are matrices
with $\dim\rho$ rows and $\dim\sigma$ columns. The coupled basis consists of
matrices having $1$ in the $i$th row and $j$th column, and all other entries
equal to $0$. Denote by $c^{k[i,j]}_{\tau[\rho,\sigma]}$ the coefficients of
expansion of the vectors of uncoupled basis in the coupled basis:
\begin{equation}  \label{eq:ClebschGordan}
\mathbf{e}^{\tau}_k=\sum_{i=1}^{\dim\rho}\sum^{\dim\sigma}_{j=1} c^{k[i,j]}_{%
\tau[\rho,\sigma]}\mathbf{e}^{\rho}_i\otimes\mathbf{e}^{\sigma}_j.
\end{equation}
The numbers $c^{k[i,j]}_{\tau[\rho,\sigma]}$ are called the \emph{%
Clebsch--Gordan coefficients} of the group~$G$. In the coupled basis, the
vectors of the uncoupled basis are matrices $c^k_{\tau[\rho,\sigma]}$ with
matrix entries $c^{k[i,j]}_{\tau[\rho,\sigma]}$, the \emph{Clebsch--Gordan
matrices}.

\begin{example}[Irreducible unitary representations of $\mathrm{SU}(2)$]
Let $\ell$ be a non-negative integer or half-integer (the half of an odd
integer) number. Let $(\rho_0,\mathbb{C}^1)$ be the trivial representation,
and let $(\rho_{1/2},\mathbb{C}^2)$ be the defining representation of $%
\mathrm{SU}(2)$. The representation $(\rho_{\ell},\mathbb{C}^{2\ell+1})$
with $\ell=1$, $3/2$, $2$, \dots, is the symmetric tensor power $\rho_{\ell}=%
\mathsf{S}^{2\ell}(\rho_{1/2})$. No other irreducible unitary
representations exist.

We may realise the representations $\rho_{\ell}$ in the space $\mathcal{P}%
^{2\ell}(\mathbb{C}^2)$ of homogeneous polynomials of degree $2\ell$ in two
formal complex variables $\xi$ and $\eta$ over the two-dimensional complex
linear space $\mathbb{C}^2$. The group $\mathrm{SU}(2)$ consists of the
matrices
\begin{equation}  \label{eq:su2}
g=
\begin{pmatrix}
\alpha & \beta \\
-\overline{\beta} & \overline{\alpha}%
\end{pmatrix}
,\qquad\alpha,\beta\in\mathbb{C},\quad|\alpha|^2+|\beta|^2=1.
\end{equation}
The representation $\rho_{\ell}$ acts as follows:
\begin{equation*}
(\rho_{\ell}(g)h)(\xi,\eta)=h(\overline{\alpha}\xi-\beta\eta, \overline{\beta%
}\xi+\alpha\eta),\qquad h\in\mathcal{P}^{2\ell}(V).
\end{equation*}
Note that $\rho_{\ell}(-E)=E$ if and only if $\ell$ is integer.

The \emph{Wigner orthonormal basis} in the space $\mathcal{P}^{2\ell}(V)$ is
as follows:
\begin{equation}  \label{eq:Wigner}
\mathbf{e}_m(\xi,\eta):=(-1)^{\ell+m}\sqrt{\frac{(2\ell+1)!} {%
(\ell+m)!(\ell-m)!}}\xi^{\ell+m}\eta^{\ell-m},\qquad m=-\ell,-\ell+1,\dots,
\ell.
\end{equation}
The matrix entries of the operators $\rho_{\ell}(g)$ in the above basis are
called \emph{Wigner $D$ functions} and are denoted by $D^{\ell}_{mn}(g)$.
The tensor product $\rho_{\ell_1}\otimes\rho_{\ell_2}$ is expanding as
follows
\begin{equation*}
\rho_{\ell_1}(g)\otimes\rho_{\ell_2}(g)
=\sum^{\ell_1+\ell_2}_{\ell=|\ell_1-\ell_2|}\oplus\rho_{\ell}(g).
\end{equation*}
\end{example}

\begin{example}[Irreducible unitary representations of $\mathrm{SO}(3)$ and $%
\mathrm{O}(3)$]
Realise the linear space $\mathbb{R}^3$ with coordinates $x_{-1}$, $x_0$,
and $x_1$ as the set of traceless Hermitian matrices over $\mathbb{C}^2$
with entries
\begin{equation*}
\begin{pmatrix}
x_0 & x_1+\mathrm{i}x_{-1} \\
x_1-\mathrm{i}x_{-1} & -x_0%
\end{pmatrix}
.
\end{equation*}
The matrix \eqref{eq:su2} acts on the so realised $\mathbb{R}^3$ as follows:
\begin{equation*}
\pi(g)
\begin{pmatrix}
x_0 & x_1+\mathrm{i}x_{-1} \\
x_1-\mathrm{i}x_{-1} & -x_0%
\end{pmatrix}
:=g^*
\begin{pmatrix}
x_0 & x_1+\mathrm{i}x_{-1} \\
x_1-\mathrm{i}x_{-1} & -x_0%
\end{pmatrix}
g.
\end{equation*}
The mapping $\pi$ is a homomorphism of $\mathrm{SU}(2)$ onto $\mathrm{SO}(3)$%
. The kernel of $\pi$ is $\pm E$. Assume that $(\rho,V)$ is an irreducible
unitary representation of $\mathrm{SO}(3)$. Then $(\rho\circ\pi,V)$ is an
irreducible unitary representation of $\mathrm{SU}(2)$ with kernel $\pm E$.
Then we have $\rho\circ\pi=\rho_{\ell}$ for some integer $\ell$. In other
words, every irreducible unitary representation $(\rho_{\ell},V)$ of $%
\mathrm{SU}(2)$ with integer $\ell$ gives rise to an irreducible unitary
representation of $\mathrm{SO}(3)$, and no other irreducible unitary
representations exist. We denote the above representation of $\mathrm{SO}(3)$
again by $(\rho_{\ell},V)$.

Let $\mathrm{SO}(2)$ be the subgroup of $\mathrm{SO}(3)$ that leaves the
vector $(0,0,1)^{\top}$ fixed. The restriction of $\rho_{\ell}$ to $\mathrm{%
SO}(2)$ is equivalent to the direct sum of irreducible unitary
representations $(\mathrm{e}^{\mathrm{i}m\varphi},\mathbb{C}^1)$, $-\ell\leq
m\leq\ell$ of $\mathrm{SO}(2)$. Moreover, the space of the representation $(%
\mathrm{e}^{\mathrm{i}m\varphi},\mathbb{C}^1)$ is spanned by the vector $%
\mathbf{e}_m(\xi,\eta)$ of the Wigner basis \eqref{eq:Wigner}.This is where
their enumeration comes from.

The group $O(3)$ is the Cartesian product of its normal subgroups $\mathrm{SO%
}(3)$ and $\{I,-I\}$. The elements of $\mathrm{SO}(3)$ are rotations,%
\index{rotation} while the elements of the second component are reflections.%
\index{reflection} Therefore, any irreducible unitary representation of $%
O(3) $ is the outer tensor product of some $(\rho_{\ell},V)$ by an
irreducible unitary representation of $\{E,-E\}$. The latter group has two
irreducible unitary representation: trivial $(\rho_+,\mathbb{C}^1)$ and
determinant $(\rho_-,\mathbb{C}^1)$. Denote $\rho_{\ell,+}:=\rho_{\ell}%
\hat{\otimes}\rho_+$ and $\rho_{\ell,-}:=\rho_{\ell}\hat{\otimes}\rho_-$.
These are all irreducible unitary representations of $O(3)$.

Introduce the coordinates on $\mathrm{SO}(3)$, the \emph{Euler angles}. Any
rotation $g$ may be performed by three successive rotations:

\begin{itemize}
\item rotation $g_0(\psi)$ about the $x_0$-axis through an angle $\psi$, $%
0\leq\psi<2\pi$;

\item rotation $g_{-1}(\theta)$ about the $x_{-1}$-axis through an angle $%
\theta$, $0\leq\theta\leq\pi$,

\item rotation $g_0(\varphi)$ about the $x_0$-axis through an angle $\varphi$%
, $0\leq\varphi<2\pi$.
\end{itemize}

The angles $\psi$, $\theta$, and $\varphi$ are the \emph{Euler angles}. The
Wigner $D$ functions are $D^{\ell}_{mn}(\varphi,\theta,\psi)$. The Wigner $D$
functions $D^{\ell}_{m0}$ do not depend on $\psi$ and may be written as $%
D^{\ell}_{m0}(\varphi,\theta)$. The \emph{spherical harmonics}%
\index{spherical harmonics} $Y_{\ell}^m$ are defined by
\begin{equation}  \label{eq:sphericalcomplex}
Y_{\ell}^m(\theta,\varphi):=%
\sqrt{\frac{2\ell+1}{4\pi}} \overline{D^{\ell}_{m0}(\varphi,\theta)}.
\end{equation}
Let $(r,\theta,\varphi)$ be the spherical coordinates in $\mathbb{R}^3$:
\begin{equation}  \label{eq:spherical}
\begin{aligned} x_{-1}&=r\sin\theta\sin\varphi,\\ x_0&=r\cos\theta,\\
x_1&=r\sin\theta\cos\varphi. \end{aligned}
\end{equation}
The measure $\mathrm{d}\Omega:=\sin\theta\,\mathrm{d}\varphi\,\mathrm{d}%
\theta$ is the Lebesgue measure on the \emph{unit sphere}%
\index{unit sphere} $S^2:=\{\,\mathbf{x}\in\mathbb{R}^3\colon\|\mathbf{x}%
\|=1\,\}$. The spherical harmonics are orthonormal:
\begin{equation*}
\int_{S^2}Y_{\ell_1}^{m_1}(\theta,\varphi)%
\overline{Y_{\ell_2}^{m_2}(\theta,\varphi)} \,\mathrm{d}\Omega=\delta_{%
\ell_1\ell_2}\delta_{m_1m_2}.
\end{equation*}
\end{example}

\begin{example}[Irreducible orthogonal representations of $\mathrm{SO}(3)$
and $\mathrm{O}(3)$]
\label{ex:o3orthogonal}

The first model is as follows. For any polynomial $h\in\mathcal{P}^{2\ell}(%
\mathbb{C}^2)$ denote by $\overline{h}$ the polynomial whose coefficients
are conjugate to those of $h$. Define the mapping $J\colon\mathcal{P}%
^{2\ell}(\mathbb{C}^2)\to\mathcal{P}^{2\ell}(\mathbb{C}^2)$ as
\begin{equation*}
(Jh)(\xi,\eta):=\overline{h}(-\eta,\xi).
\end{equation*}
The orthonormal basis of eigenvectors of $J$ with eigenvalue~$1$ was
proposed by \citet{MR1888117}. The vectors of the \emph{Gordienko basis} are
as follows ($m\geq 1$):
\begin{equation*}
\begin{aligned}
\mathbf{h}_{-m}(\xi,\eta)&:=\frac{(-\mathrm{i})^{\ell-1}}{\sqrt{2}}
[(-1)^m\mathbf{e}_m(\xi,\eta)-\mathbf{e}_{-m}(\xi,\eta)],\\
\mathbf{h}_0(\xi,\eta)&:=(-\mathrm{i})^{\ell}\mathbf{e}_0(\xi,\eta),\\
\mathbf{h}_m(\xi,\eta)&:=-\frac{(-\mathrm{i})^{\ell}}{\sqrt{2}}
[(-1)^m\mathbf{e}_m(\xi,\eta)+\mathbf{e}_{-m}(\xi,\eta)]. \end{aligned}
\end{equation*}
In this basis, the representations $\rho_{\ell,+}$ and $\rho_{\ell,-}$
become orthogonal and will be denoted by $U^{\ell g}$ and $U^{\ell u}$ (g by
German \emph{gerade},%
\index{representation!gerade} even, and u by \emph{ungerade},%
\index{representation!ungerade} odd).

The Clebsh--Gordan coefficients of the groups $\mathrm{SO}(3)$ and $\mathrm{O%
}(3)$ with respect to the Gordienko basis were calculated by %
\citet{MR2078714}. We call them \emph{Godunov--Gordienko coefficients} and
denote them by $g^{m[m_1,m_2]}_{%
\ell[\ell_1,\ell_2]}$. An algorithm for calculation of the
Godunov--Gordienko coefficients was proposed by \citet{MR3308053}.
\end{example}

\begin{example}[Expansions of tensor representations of the group $\mathrm{O}%
(3)$]
\label{ex:8}

Let $r\geq 2$ be a nonnegative integer, and let $\Sigma_r$ be the
permutation group of the numbers $1$, $2$, \dots, $r$. The action
\begin{equation*}
\sigma\cdot(\mathbf{v}_1\otimes\cdots\otimes\mathbf{v}_r):=\mathbf{v}%
_{\sigma^{-1}(1)} \otimes\cdots\otimes\mathbf{v}_{\sigma^{-1}(r)},\qquad%
\sigma\in\Sigma_r,
\end{equation*}
may be extended by linearity to an orthogonal representation of the group $%
\Sigma_r$ in the space $(\mathbb{R}^3)^{\otimes r}$, call it $(\rho_r,(%
\mathbb{R}^3)^{\otimes r})$. Consider the orthogonal representation $(\tau,(%
\mathbb{R}^3)^{\otimes r})$ of the group $\mathrm{O}(3)\times\Sigma_r$
acting by
\begin{equation*}
\tau(g,\sigma)(\mathsf{T}):=\rho^{\otimes r}(g)\rho_r(\sigma)(\mathsf{T}%
),\qquad \mathsf{T}\in(\mathbb{R}^3)^{\otimes r}.
\end{equation*}

The representation $(\tau,(\mathbb{R}^3)^{\otimes 2})$ of the group $\mathrm{%
O}(3)\times\Sigma_2$ is the direct sum of three irreducible components
\begin{equation*}
\tau=[\rho_{0,+}(g)\tau_+(\sigma)]\oplus[\rho_{1,+}(g)\varepsilon(\sigma)]
\oplus[\rho_{2,+}(g)\tau_+(\sigma)],
\end{equation*}
where $\tau_+$ is the trivial representation of the group $\Sigma_2$, while $%
\varepsilon$ is its non-trivial representation. The one-dimensional space of
the first component is the span of the identity matrix and consists of
scalars. The three-dimensional space of the second component is the space $%
\mathsf{\Lambda}^2(\mathbb{R}^3)$ of $3\times 3$ skew-symmetric matrices.
Its elements are three-dimensional pseudo-vectors. Finally, the
five-dimensional space of the third component consists of $3\times 3$
traceless symmetric matrices (deviators). The second component is $(\mathsf{%
\Lambda}^2(g),\mathsf{\Lambda}^2(\mathbb{R}^3))$, and the direct sum of the
first and third components is $(\mathsf{S}^2(g),\mathsf{S}^2(\mathbb{R}^3))$.

In general, the representation $(\tau,(\mathbb{R}^3)^{\otimes r})$ is
reducible and may be represented as the direct sum of irreducible
representations as follows:
\begin{equation*}
\tau(g,\sigma)=\sum_{\ell=0}^{r}\sum_{q=1}^{N^{\ell}_r}\oplus U^{\ell x}(g)
\rho_q(\sigma),
\end{equation*}
where $q$ is called the \emph{seniority index} of the component $U^{\ell
x}(g)\rho_q(\sigma)$, see \citet{PhysRevA.25.2647}, and where $x=g$ for even
$r$ and $x=u$ for odd $r$. The number $N^{\ell}_r$ of copies of the
representation $U^{\ell x}$ is given by
\begin{equation*}
N^{\ell}_r=\sum_{k=0}^{\llcorner(r-\ell)/3\lrcorner}\binom{r}{k} \binom{%
2r-3k-\ell-2}{r-2}.
\end{equation*}
\end{example}

\section{Classical invariant theory}

\label{ap:invariant}

Let $V$ and $W$ be two finite-dimensional linear spaces over the same field $%
\mathbb{K}$. Let $(\rho,V)$ and $(\sigma,W)$ be two representations of a
group~$G$. A mapping $h\colon W\to V$ is called a \emph{covariant} or \emph{%
form-invariant} or a \emph{covariant tensor} of the pair of representations $%
(\rho,V)$ and $(\sigma,W)$, if
\begin{equation*}
h(\sigma(g)\mathbf{w})=\rho(g)h(\mathbf{w}),\qquad g\in G.
\end{equation*}
In other words, the diagram
\begin{equation*}
\xymatrix{ W\ar[r]^h\ar[d]^{\sigma} & V\ar[d]^{\rho} \\ W\ar[r]^h & V }
\end{equation*}
is commutative.

If $V=\mathbb{F}^1$ and $\rho$ is the trivial representation of $G$, then
the corresponding covariant scalars are called \emph{absolute invariants}
(or just invariants) of the representation $(\sigma,W)$, hence the name
\emph{Invariant Theory}. Note that the set $\mathbb{K}[W]^G$ of invariants
is an \emph{algebra} over the field $\mathbb{K}$, that is, a linear space
over $\mathbb{F}$ with bilinear multiplication operation and a
multiplication identity $1$. The product of a covariant $h\colon W\to V$ and
an invariant $f\in\mathbb{K}[W]^G$ is again a covariant. In other words, the
covariant tensors of the pair of representations $(\rho,V)$ and $(\sigma,W)$
form a \emph{module} over the algebra of invariants of the representation $%
(\sigma,W)$.

A mapping $h\colon W\to V$ is called \emph{homogeneous polynomial mapping of
degree~$d$} if for any $\mathbf{v}\in V$ the mapping $\mathbf{w}\mapsto (h(%
\mathbf{w}),\mathbf{v})$ is a homogeneous polynomial of degree~$d$ in $\dim
W $ variables. The mapping $h$ is called a \emph{polynomial covariant of
degree $d$} if it is homogeneous polynomial mapping of degree~$d$ and a
covariant.

Let $(\sigma,W)$ be the defining representation of $G$, and $(\rho,V)$ be
the $r$th tensor power of the defining representation. The corresponding
covariant tensors are said to have \emph{an order}~$r$. The covariant
tensors of degree~$0$ and of order~$r$ of the group $\mathrm{O}(n)$ are
known as \emph{isotropic tensors}.

The algebra of invariants and the module of covariant tensors were an object
of intensive research. The first general result was obtained by %
\citet{Gordan1868}. He proved that for any finite-dimensional complex
representation of the group $G=\mathrm{SL}(2,\mathbb{C})$ the algebra of
invariants and the module of covariant tensors are finitely generated. In
other words, there exists an \emph{integrity basis}: a finite set of
invariant homogeneous polynomials $I_1$, \dots, $I_N$ such that every
polynomial invariant can be written as a polynomial in $I_1$, \dots, $I_N$.
An integrity basis is called \emph{minimal} if none of its elements can be
expressed as a polynomial in the others. A minimal integrity basis is not
necessarily unique, but all minimal integrity bases have the same amount of
elements of each degree.

The algebra of invariants is not necessarily free. Some polynomial relations
between generators, called \emph{syzygies} may exist.

The importance of polynomial invariants can be explained by the following
result. Let $G$ be a closed subgroup of the group $\mathrm{O}(3)$, the group
of symmetries of a material. Let $(\rho,\mathsf{V})$, $(\rho_1,\mathsf{V}_1)$%
, \dots, $(\rho_N,\mathsf{V}_N)$ be finitely many orthogonal representations
of $G$ in real finite-dimensional spaces. Let $\mathsf{T}\colon\mathsf{V}%
_1\oplus\cdots\oplus\mathsf{V}_N\to\mathsf{V}$ be an \emph{arbitrary} (say,
measurable) covariant of the pair $\rho$ and $\rho_1\oplus\cdots\oplus\rho_N$%
. Let $\{\,I_k\colon 1\leq k\leq K\,\}$ be an integrity basis for \emph{%
polynomial} invariants of the representation $\rho$, and let $\{\,\mathsf{T}%
_l\colon 1\leq l\leq L\,\}$ be an integrity basis for \emph{polynomial}
covariant tensors of the pair $\rho$ and $\rho_1\oplus\cdots\oplus\rho_N$.
Following \citet{MR0171421}, we call $\mathsf{T}_l$ \emph{basic covariant
tensors}.

\begin{theorem}[\citet{MR0171421}]
\label{th:Wineman-Pipkin} A function $\mathsf{T}\colon\mathsf{V}%
_1\oplus\cdots\oplus\mathsf{V}_N\to\mathsf{V}$ is a measurable covariant of
the pair $\rho$ and $\rho_1\oplus\cdots\oplus\rho_N$ if and only if it has
the form
\begin{equation*}
\mathsf{T}(\mathsf{T}_1,\dots,\mathsf{T}_N)=\sum_{l=1}^{L}\varphi_l(I_1,%
\dots,I_K) \mathsf{T}_l(\mathsf{T}_1,\dots,\mathsf{T}_N),
\end{equation*}
where $\varphi_l$ are real-valued measurable functions of the elements of an
integrity basis.
\end{theorem}

In 1939 in the first edition of \citet{MR1488158} Hermann Weyl proved that
any polynomial covariant of degree~$d$ and of order~$r$ of the group $%
\mathrm{O}(n)$ is a linear combination of products of Kronecker's deltas $%
\delta_{ij}$ and second degree homogeneous polynomials $x_ix_j$.

\end{document}